\documentclass[11pt]{article}
\usepackage{amsmath, amsfonts, amssymb, amscd, enumerate, eucal, tikz, stmaryrd}
\usepackage{color}
\usepackage{graphicx}
\usepackage{verbatim,amscd}
\usepackage[british]{babel}
\usepackage{amstext, float}
\usepackage[utf8]{inputenc}
\usepackage{array,color}
\usepackage{enumitem}
\usepackage[small]{eulervm}
\usepackage[unicode]{hyperref}
\usepackage{pgf,tikz}
\usepackage{mathrsfs}
\usetikzlibrary{arrows}
\definecolor{xdxdff}{rgb}{0.49019607843137253,0.49019607843137253,1.}
\definecolor{qqqqff}{rgb}{0.,0.,1.}
\definecolor{cqcqcq}{rgb}{0.7529411764705882,0.7529411764705882,0.7529411764705882}
\definecolor{ttqqqq}{rgb}{0.2,0.,0.}
\def\a{\alpha}

\def\i{\iota}

\def\k{\kappa}

\def\o{\omega}
\def\O{\Omega}

\def\u{\upsilon}

\def\O{\Omega}

\chardef\tempcat=\the\catcode`\@ \catcode`\@=11
\def\cyracc{\def\u##1{\if \i##1\accent"24 i
    \else \accent"24 ##1\fi }}
\newfam\cyrfam

\DeclareFontFamily{OT1}{msb}{}{} \DeclareFontShape{OT1}{msb}{m}{n}
 {  <5> <6> <7> <8> <9> <10> gen * msbm
      <10.95><12><14.4><17.28><20.74><24.88>msbm10}{}
\DeclareMathAlphabet{\bubble}{OT1}{msb}{m}{n}
\def\bR{{\mathbb R}}
\def\bZ{{\mathbb Z}}
\def\bC{{\mathbb C}}

\newcommand{\Be}{\boldsymbol{e}}
\newfont{\goth}{eufm10 scaled \magstep1}
\def\ga{\mathfrak a}
\def\gb{\mathfrak b}

\def\gf{\mathfrak f}
\def\gg{\mathfrak g}

\def\gh{\mathfrak h}
\def\gk{\mathfrak k}

\def\gm{\mathfrak m}
\def\gn{\mathfrak n}

\def\gp{\mathfrak p}

\def\gs{\mathfrak s}

\def\gu{\mathfrak u}

\def\gso{\mathfrak{so}}

\def\gsp{\mathfrak{sp}}

\newcommand{\fn}{\mathfrak n}
\newcommand{\fb}{\mathfrak b}
\newcommand{\ff}{\mathfrak f}
\newcommand{\rk}{\operatorname{rk}}
\newcommand{\stab}{\mathfrak{stab}}

\newcommand{\ggl}{\mathfrak{gl}}
\newcommand{\fsp}{\mathfrak{sp}}
\newcommand{\fsl}{\mathfrak{sl}}
\def\gsl{\mathfrak{sl}}
\newcommand{\fso}{\mathfrak{so}}

\newcommand{\fp}{\mathfrak{p}}
\newcommand{\fs}{\mathfrak{s}}
\newcommand{\fu}{\mathfrak{u}}
\newcommand{\fh}{\mathfrak{h}}
\newcommand{\fk}{\mathfrak{k}}

\newcommand{\RR}{\mathbb{R}}

\newfont{\mcal}{eusm10 scaled \magstep1}

\def\p{\partial}

\def\exp{\mathrm{exp}}

\def\Id{\mathrm{Id}}
\def\tr{\mathrm{tr\;}}

\def\Ad{\mathrm{Ad}}
\def\ad{\mathrm{ad}}
\def\diag{\mathrm{diag \;}}

\def\exp{ \mathrm exp}
\def\span{ \mathrm span \;}
\newtheorem{Th}{Theorem}
\newtheorem{Prop}{Proposition}
\newtheorem{Cor}{Corollary}
\newtheorem{Lem}{Lemma}
\newtheorem{Def}{Definition}
\newtheorem{Rem}{Remark}
\def\bt{\begin{Th}}
\def\et{\end{Th}}
\def\bp{\begin{Prop}}
\def\ep{\end{Prop}}
\def\bc{\begin{Cor}}
\def\ec{\end{Cor}}
\def\bl{\begin{Lem}}
\def\el{\end{Lem}}
\def\bd{\begin{Def}}
\def\ed{\end{Def}}
\def\br{\begin{Rem}}
\def\er{\end{Rem}}
\def\pf{\noindent{\it Proof. }}
\def\qed{\hspace{2ex} \hfill $\square $ \par \medskip}

\def\be{\begin{equation}}
\def\ee{\end{equation}}
\def\arr{\begin{array}{rlll}}
\def\ea{\end{array}}
\def\bea{\begin{eqnarray}}
\def\eea{\end{eqnarray}}
\def\bean{\begin{eqnarray*}}
\def\eean{\end{eqnarray*}}

\def\span {\mathrm {span} \, }

\def\diag{ \mathrm {diag} \,}

\begin{document}
\title{Homogeneous symplectic 4-manifolds  and finite \\ dimensional  Lie  algebras of  symplectic vector  fields \\ on  the symplectic  4-space}
\author{D.\ V.\ Alekseevsky, A.\ Santi}
\date{}
\maketitle
\begin{abstract}
We classify the finite type (in the sense of E.\ Cartan theory of prolongations) subalgebras $\mathfrak h\subset\fsp(V)$, where $V$ is the symplectic $4$-dimensional space, and show that they satisfy $\mathfrak h^{(k)}=0$ for all $k>0$. Using this result, we
reduce the problem of classification of graded transitive finite-dimensional Lie algebras $\mathfrak g$ of symplectic vector fields on $V$ to the description of graded transitive finite-dimensional subalgebras of the full prolongations $\gp_1^{(\infty)}$ and $\gp_2^{(\infty)}$, where $\gp_1$ and $\gp_2$ are the maximal parabolic subalgebras of $\gsp(V)$. We then classify all such $\gg\subset\gp_i^{(\infty)}$, $i=1,2$, under some assumptions, and describe the associated $4$-dimensional homogeneous symplectic manifolds $(M = G/K, \omega)$. We prove that any reductive homogeneous symplectic manifold (of any dimension) admits an invariant torsion free symplectic connection, i.e., it is a homogeneous Fedosov manifold, and give conditions for the uniqueness of the Fedosov structure. Finally, we show that any nilpotent symplectic Lie group (of any dimension) admits 
a natural invariant Fedosov structure which is Ricci-flat.
\end{abstract}
\tableofcontents
\section{Introduction, statement of   the    results}
Let   $(M = G/K, \omega)$ be   a   homogeneous  symplectic  manifold   of  a  connected  Lie  group $G$.
 If  the   group  $G$ is  fixed,  the  description of  all    homogeneous  $G$-manifolds   reduces  to  the description  of the closed  2-forms
 $\Omega \in  \Lambda^2_{cl}\mathfrak{g}^*$ on  its Lie   algebra. In this paper, we  consider  another problem:  to describe  the  homogeneous  symplectic manifolds of  a  given  dimension $\dim M=2n$.

	We will consider  the  following  approach. The  Lie    algebra  $\mathfrak{g} = \mathrm{Lie}(G)$,  considered  as a  finite  dimensional Lie   algebra of symplectic vector  fields  on the symplectic manifold $(M=G/K, \omega)$,  admits  a  natural  filtration
\begin{equation}
\label{eq:filteredg}
\mathfrak{g}= \mathfrak{g}_{-1} \supset \mathfrak{g}_0 \supset \mathfrak{g}_1 \supset \cdots \supset \mathfrak{g}_k \supset \{0  \}\;,
\end{equation}
	given by the symplectic vector fields vanishing at a certain order at $o=eK$.
 We consider the  associated  graded  Lie  algebra
  \begin{equation}
	\label{eq:graded}
	gr (\mathfrak{g}) = \mathfrak{g}^{-1} + \mathfrak{g}^0 + \mathfrak{g}^1 + \cdots + \mathfrak{g}^k\;,
	\end{equation}
  it is   a  transitive  graded subalgebra   of  the  Lie   algebra  $\mathfrak{sp}(V)^{(\infty)}  =   V + \mathfrak{sp}(V) + \mathfrak{sp}(V)^{(1)} + \cdots$
  of all formal  symplectic vector   fields, where  $V = \mathfrak{g}^{-1} = T_oM,\,  o = eK$.
  We note   that the Lie algebra $\gg$ can be realized as a filtered deformation of the graded Lie algebra $gr(\gg)$ and that  $\mathfrak{g}^j$ is  a  subspace  of  the $j$-th  prolongation   $\mathfrak{h}^{(j)}$  of  the  linear   isotropy   subalgebra   $\mathfrak{h} = \mathfrak{g}^0  \subset  \mathfrak{sp}(V)$. 
	
  The    first    arising  question is  to  describe   all   subalgebras   $\mathfrak{h} \subset \mathfrak{sp}(V)$  of  finite  type,  i.e.,    algebras  which have  finite-dimensional  full prolongation
  $$  \mathfrak{h}^{(\infty)} = V + \mathfrak{h} + \mathfrak{h}^{(1)} + \cdots + \mathfrak{h}^{(k)}\;.$$
For  simplicity, consider at first   the  complex  case,  where  $\mathfrak{h}$  is  a   complex  subalgebra  of   the   complex  symplectic   Lie algebra  of   the    symplectic vector  space   $V = \mathbb{C}^{2n}$.
  Then    a   complex  linear  Lie algebra   $\mathfrak{h}$  is of  finite    type if and only if it  has  no  rank one   endomorphisms.
  For $n=1$, $\mathfrak{sp}(V) = \mathfrak{sl}_2(\mathbb{C})$  and   (up  to a  conjugation)   the  Cartan   subalgebra  $\mathbb{C} \mathrm{diag}(1,-1)$  is  the only  nonzero  subalgebra of   finite  type. In this paper, we  deal with  the   case  $n=2$, paving the way to the classification of homogeneous symplectic $4$-manifolds and homogeneous Fedosov $4$-manifolds.

In \S\ref{sec:maximal} and \S\ref{sec:classification}, we describe the maximal subalgebras of $\fsp_2(\bR)$ and determine those of finite type. We then reduce the classification problem of (maximal) finite type subalgebras of $\fsp_2(\bR)$ to the description of finite type subalgebras of maximal parabolic subalgebras $\fp_1$, $\fp_2$ and the subalgebra $\fs_1=\fsp(V_1)\oplus\fsp(V_2)$ preserving an orthogonal decomposition $V=V_1+V_2$ of $V$. Theorem \ref{thm:1} below describes the finite type subalgebras of $\fsp_2(\bR)$ and it is a consequence of this analysis, which is carried out in \S\ref{sec:classification}.

Throughout the paper, we shall identify $\fsp_2(\bR)$  with
the space of quadratic polynomials on $V^*\cong V=\mathrm{span}(p_1,p_2,q_1,q_2)$.
\bt
\label{thm:1}
Let $\gh$ be a subalgebra of $\fsp_2(\bR)$. Then $\gh$ is either of infinite type or it is, up to conjugation, a subalgebra of one of the Lie algebras of the following list:
\begin{itemize}
\item[1.] the unitary algebra $\fu_2$ (maximal compact subalgebra of $\fsp_2(\bR)$);
\item[2.] the pseudo-unitary algebra $\fu_{1,1}=  \fsl(H)\oplus  \fso(E)   \subset \fsp(H \otimes E)$,
where   the    symplectic   structure on  $V = H \otimes E = \bR^2 \otimes \bR^2$  is   the  tensor product $\Omega = {\mathrm{vol}_H} \otimes \eta_E$ of      the volume  form  on  $H$  and       an  Euclidean metric on  $E$;
\item[3.] the  reductive  subalgebra $Q\vee P\cong\ggl(P)$, where $P$ and $Q$ are  complementary  Lagrangian subspaces;
\item[4.] the irreducible singular subalgebra $\fsl_2^4(\bR)$ acting on $ V =  S^3(\bR^2)$;
\item[5.] the direct sum $\fso_2(\bR)\oplus\bR\mathrm{diag}(1,-1)$ of two Cartan subalgebras of $\fsl_2(\bR)$, one compact and one noncompact; 
\item[6.] the solvable $2$-dimensional subalgebra
\begin{equation*}
\begin{aligned}
D_{4,12}&=\mathrm{span}(p_1q_1+\epsilon p_2^2,p_2q_1)\;,
\end{aligned}
\end{equation*}
where $\epsilon=\pm 1$;
\item[7.] the  $1$-dimensional subalgebra $\mathrm{span}(p_2^2+q_2^2+\epsilon p_1^2)$, where $\epsilon=\pm 1$.
\end{itemize}
All   these   subalgebras  are  maximal     finite   type   subalgebras. In particular a finite type subalgebra of $\fsp_2(\bR)$ has at most dimension $4$.   
\et
For more details, we refer the reader to Theorems \ref{thm:3}--\ref{thm:5} and the proof in \S\ref{sec:finalthm1} of the maximality of the finite type subalgebras of Theorem \ref{thm:1}. The  upper   index in $\gsl_2^4(\bR)$ indicates  the  dimension of the irreducible  representation; the notation for the $2$-dimensional solvable subalgebra $DF_{4,12}$ makes contact with the description of conjugacy classes of subalgebras of the similitude algebra in \cite{PWSZ}.

A direct but rather surprising corollary of Theorem \ref{thm:1} is that all finite type subalgebras of $\fsp_2(\bR)$ have trivial prolongations.
\bc
\label{thm:2}
Let $\gh$ be a finite type subalgebra of $\fsp_2(\bR)$. Then the first prolongation $\fh^{(1)}=0$.
\ec
Due to this, the finite type subalgebras are a good class of candidates to construct homogeneous Fedosov manifolds. We will expand on this later.
\vskip0.2cm
We thank B.\ Kruglikov, who informed us that there is no upper bound on the dimension of a Lie group $G$ acting transitively on a symplectic $4$-manifold. This is in sharp contrast with the Riemannian case. Let $(M=G/K,\omega)$ be a homogeneous symplectic $4$-manifold with associated graded Lie algebra \eqref{eq:graded}. It is a finite-dimensional  subalgebra  of  the  full prolongation $\gh^{(\infty)}$ of 
the  linear   isotropy   subalgebra   $\mathfrak{h} = \mathfrak{g}^0  \subset  \mathfrak{sp}(V)$. If $\gh$ has finite type then $\gh^{(1)}=0$ by Corollary \ref{thm:2} and therefore $gr(\gg)  = V +\gh$. If $\gh=\gsp(V)$ then the first prolongation $\gh^{(1)}$ is irreducible and, together with $V+\gh$, it generates the (infinite-dimensional) full prolongation $\mathfrak{sp}(V)^{(\infty)}$. Hence $gr(\gg)  = V +\gh$ also in this case. According to the description of maximal subalgebras (Theorem \ref{thm:maxsub} and Corollary \ref{cor:finitetype}), the  only maximal infinite type  subalgebras  of    $\gsp(V)$  are  the two   maximal parabolic  subalgebras   $\gp_1$ and $\gp_2$  and the semisimple   subalgebras  $\fs_1=\fsp(V_1)\oplus\fsp(V_2)$ and $\gs_4=\gsl_2(\bC)$.

The classification  of  all   graded  transitive  finite-dimensional   Lie  algebras  of  symplectic vector  fields  on $V$   reduces  in this way to   the description  of  the  graded transitive   finite-dimensional  subalgebras     of  the  infinite-dimensional  Lie  algebras  $\gp_1^{(\infty)}, \gp_2^{(\infty)},   \gs_1^{(\infty)}, \gs_4^{(\infty)}$.  Furthermore, using Proposition \ref{prop:S4} and the results of \S\ref{sec:classification3.1}, one  can directly see   that  such  subalgebras   of   $ \gs_1^{(\infty)}$ and $\gs_4^{(\infty)} $  are  either of the form $gr(\gg)  = V +\gh$ as above or they are subalgebras  of $\gp_1^{(\infty)}$ or $\gp_2^{(\infty)}$, see \S \ref{sec:dichotomy} for full details. In other words, we have the following dichotomy result of Theorem \ref{thm:dichotomy}. We recall that a filtered deformation of a graded Lie algebra is called trivial if it is itself graded.

\bt
\label{thm:dichotomy}
Let $(M = G/K, \omega)$ be a $4$-dimensional homogeneous symplectic manifold,  on which a finite-dimensional connected Lie group $G$ (not
necessarily compact) acts transitively and almost effectively with connected stabilizer $K\subset G$. Then the Lie algebra $\mathfrak g=Lie(G)$ is a (possibly trivial) filtered deformation of $gr(\gg)$ and
exactly one of the following two cases occurs:
\begin{enumerate}
\item $gr(\mathfrak{g})=V+\mathfrak{h}$ is the affine Lie algebra, where $\mathfrak h$ is any subalgebra of $\mathfrak{sp}(V)$; 
\item $gr(\mathfrak g)=V+\gh+\gg^1+\cdots+\gg^k$ is such that $\mathfrak{g}^1\neq 0$, where $\mathfrak h$ is an infinite type subalgebra of the maximal parabolic subalgebra $\gp_1$ or the maximal parabolic subalgebra $\gp_2$.
\end{enumerate}
If the isotropy representation $j:K\to j(K)=H\subset \mathrm{Sp}_2(\bR)$ is not infinitesimally exact then there exists an isotropic distribution $\mathcal D\subset TM$ on $M=G/K$ of $\rk\mathcal D=1$ or $2$ left invariant by $G$. 
\et

The first alternative of Theorem \ref{thm:dichotomy} parallels the situation for Riemannian manifolds. In \S \ref{sec:K} we focus on the second alternative and use a more suggestive description of the full prolongations $\gp_1^{(\infty)}$ and $\gp_2^{(\infty)}$ in order to describe their finite-dimensional subalgebras. Recall that a Lie algebra of vector fields is called primitive if it does not preserve any integrable subdistribution. The notion of a transversally primitive transitive Lie algebra of vector fields is given in \S\ref{sec:preldis}.
Under some assumptions, we describe all finite-dimensional transitive and transversally primitive subalgebras of $\gp_1^{(\infty)}$ and $\gp_2^{(\infty)}$. They are given in Theorem \ref{thm:K1} in \S \ref{sec:4.2} and in Theorem \ref{thm:K2} in \S \ref{sec:4.3}. We then construct homogeneous symplectic $4$-manifolds for any of these Lie algebras and show that most of them do not admit any invariant torsion free symplectic connection.
\vskip0.2cm\par
A Fedosov structure on a symplectic manifold $(M,\o)$ is the assignment of a torsion free symplectic connection $\nabla$ \cite{GRS}.
The celebrated result of Fedosov gives a deformation quantization canonically defined by the data $(M,\o,\nabla)$ \cite{Fed}.
Torsion free symplectic connections exist on any symplectic manifold, in other words, a structure of a symplectic manifold can always be extended to a structure of a Fedosov manifold.

The are certain classes of symplectic manifolds for which there is a natural choice of a {\it unique} Fedosov structure, for instance those preserving some additional geometric data. These include Levi-Civita connections on K\"ahler manifolds, symplectic manifolds endowed with a Lagrangian polarization \cite{Hs} and symmetric symplectic spaces \cite{Bie}.  We also note that there is a notion of a {\it preferred} Fedosov structure, which is not imposed by the presence of extra geometric structures (see \cite{BCGRS}). Every compact coadjoint orbit has an invariant preferred Fedosov structure but the connection is not unique in general, see \cite{CGR}. Two important classes of symplectic manifolds with preferred Fedosov structures are the class of {\it Ricci-type} manifolds --- the Fedosov manifolds $(M,\o,\nabla)$ for which the curvature of $\nabla$ is entirely determined by the Ricci tensor $\operatorname{ric}^{\nabla}$ --- and the {\it Ricci flat} manifolds, i.e., those satisfying $\operatorname{ric}^{\nabla}=0$. The Ricci-type connections are well-understood, both from a local and a global point of view \cite{BCGRS}. Less is known on Ricci-flat Fedosov manifolds, but there is a procedure to construct examples in dimension $\geq 6$ \cite{BCGRS}.

In \S\ref{sec:lastsecI}, we show that any {\it reductive}  homogeneous symplectic manifold $(M = G/K,\o)$ admits an invariant torsion free symplectic connection, in other words, it is a homogeneous Fedosov manifold (see Proposition \ref{prop:redFedosov}).
Corollary \ref{thm:2} then suggests a way to select a class of unique Fedosov structures on reductive homogeneous manifolds.
Let $H$ be a   subgroup  of $Sp_2(\bR)$ of  finite  type   (that is, $\gh \subset \gsp_2(\bR) $ is included in one of  the  subalgebras from  Theorem \ref{thm:1}).
 We consider the reductive  homogeneous    symplectic   $4$-manifolds  $(M = G/K,\omega)$   with  isotropy   representation    $j(K) =H$ and look for invariant torsion free (symplectic) connections compatible with the canonical reduction $\pi:P=G\to M=G/K$ of the bundle of all symplectic frames on $M$.

Combining Proposition \ref{prop:redFedosov} with Theorem \ref{thm:last} in \S\ref{sec:lastsecII}, we arrive at the following.
\bt
\label{thm:redFed}
Let $(M = G/K,\o)$ be a $4$-dimensional homogeneous symplectic manifold
with a finite type isotropy $j(K)=H\subset \mathrm{Sp}_2(\bR)$. Then the isotropy
representation $j:K\to\mathrm{Sp}_2(\bR)$ is infinitesimally exact, i.e., it has a discrete kernel and
\begin{enumerate}
\item If $M=G/K$ is reductive, there exists an invariant torsion free symplectic connection;
\item If there exists
an invariant torsion free (symplectic) connection which is compatible with the canonical reduction $\pi:P\to M$,
then    $M=G/K$ is reductive. Furthermore, such a connection is unique.
\end{enumerate}
If an invariant torsion free symplectic connection exists then the isotropy representation is exact.
\et
More generally, any torsion free (symplectic) connection on a symplectic $4$-dimensional manifold $(M,\o)$ compatible with a symplectic $H$-structure $\pi:P\to M$ whose structure group $H\subset\mathrm Sp(V)$ has finite type is unique. In the case of reductive homogeneous symplectic $4$-manifolds, there is a natural choice of $\pi:P\to M$ and $2.$ of Theorem \ref{thm:redFed} can be rephrased by saying that the associated deformation quantization (if it exists) is defined solely in terms of $(M,\o)$.
	
	Torsion free connections on $4$-manifolds whose holonomies exactly realize the representation of $H=\mathrm{SL}_2(\bR)$ on $S^3(\bR^2)$ have been considered in great detail in \cite{Bry} and homogeneous Fedosov manifolds (of any dimension) with special symplectic holonomy have been constructed  in \cite{Sch}. The classification of reductive homogeneous symplectic $4$-manifolds  $(M = G/K,\omega)$  such that the isotropy subgroup is  a (maximal) finite type subgroup of $Sp_2(\bR)$ --- and of their invariant Fedosov structures --- will be the content of future work \cite{AlSa}.
	
In \S\ref{sec:lastsecIII}, we consider the case of {\it symplectic Lie groups} (of any dimension), that is, the homogeneous symplectic manifolds $(M=G/K,\o)$ with trivial stabilizer $K=\{1\}$. Using the results of \S \ref{sec:lastsecI}-\ref{sec:lastsecII} and the theory of left-symmetric algebras, we derive the following.
\bt
Any nilpotent symplectic Lie group (of any dimension) admits an invariant torsion free symplectic connection which is Ricci flat.
\et	
See Theorem \ref{thm:nilpoliegroup} for more details. We recall that the class of nilpotent symplectic Lie groups is particularly large and that they can all be described in terms of double extensions (see \cite{MR}). The explicit classification, up to local equivalence, in dimension $\leq 6$ can be found in \cite{Kr, KGM}. We remark that structural properties of solvable and nilpotent symplectic Lie groups, as well as cotangent symplectic Lie groups, have been investigated in detail in \cite{BC}. In particular, the theory of symplectic reduction w.r.t. isotropic normal subgroups is developed and irreducible symplectic Lie groups are classified. 
\vskip0.4cm\par\noindent
{\it Acknowledgements.} We would like to thank the referee for useful comments and suggestions.
\section{Maximal  subalgebras  of $\gsp_2(\bR)$}
\label{sec:maximal}
\subsection{List of  maximal  subalgebras  of $\gsp_2(\bR) \cong  \gso(2,3)$}
  Let $(V, \O)$ be the $4$-dimensional symplectic (real) vector  space.  The  $\gsp(V)$-module  $W = \Lambda_0^2(V)$  has an invariant inner product $\eta= \Omega \wedge \Omega$  of  signature
 $(2,3)$,  which gives  an isomorphism    $  \gsp_2(\bR) = \gsp(V)  \cong  \gso(W) = \gso(2,3)$.
\\

According  to \cite{GOVIII}, the
maximal   subalgebras  of the complexification  $\gso_5(\bC)$   are  exhausted  by
\begin{enumerate}
\item maximal  parabolic subalgebras,
\item the reducible subalgebra  $\gso_{4}(\bC)$,
\item a simple irreducible subalgebra  $\gsl_2^5(\bC)$
    defined  by  the   natural  action on   $S^4(\bC^2)$.
		\end{enumerate}
		We remark that the reducible subalgebra $\gso_3(\bC) \oplus \gso_{2}(\bC)$ has to be removed from the list of \cite[Theorem 3.1, p. 205]{GOVIII}, since it is contained in the parabolic subalgebra stabilizing an isotropic line, hence it is not maximal. The  upper   index in $\gsl_2^5(\bC)$ indicates  the  dimension of the irreducible  representation.

Since  $\mathfrak{so}(2,3)$  is  the   normal  real  form of  $\mathfrak{so}_5(\mathbb{C})$, there is  no  big  difference between maximal  real subalgebras of  $\mathfrak{so}(2,3)$  and  complex  subalgebras of  $\mathfrak{so}_5(\mathbb{C})$. It turns out that they are
real forms of the maximal complex subalgebras listed above or real forms of $\gso_3(\bC) \oplus \gso_{2}(\bC)$ which do not preserve any isotropic line (see \cite{PSWZ}). 

Besides maximal parabolic subalgebras, the real  forms of  maximal   complex   subalgebras  of  $\gso_5(\bC)$ which  belong  to   $\mathfrak{\gso}(2,3)\cong \fsp_2(\bR)$  give  the  following list  of  maximal  subalgebras:
 \begin{itemize}
\item[(i)] $\gso(2,2)\cong \gsl_2(\bR)\oplus\gsl_2(\bR)=\gsp(V_1) \oplus \gsp(V_2)$,
\item[(ii)] $\gso(1,3)\cong\gsl_2(\bC)$,
\item[(iii)] $\gsl_2^5(\bR)$.
\end{itemize}
In the first case, the symplectic space has an $\O$-orthogonal decomposition $V=V_1+V_2$, where $V_i = \mathrm{span}(p_i,q_i)$ for $i=1,2$.
We denote the associated semisimple  subalgebra by $\gs_1:=\gsp(V_1)\oplus\gsp(V_2)$ and note that $\gsp_2(\bR)$ has a symmetric  decomposition $\gsp_2(\bR) = \gs_1+ V_1 \vee V_2$ with corresponding pseudo-Riemannian irreducible   symmetric  space
$\mathrm{Sp}_2(\bR)/ \mathrm{SL}_2(\bR)\cdot\mathrm{SL}_2(\bR)$.

The spin representation of $\gs_4:=\fso(1,3)\cong\gsl_2(\bC)$ is symplectic real $4$-dimensional, and it consists of all symplectic matrices compatible with a complex structure $J$ such that $J^*\O=-\O$.

The last case gives the singular subalgebra $\gs_5:=\gsl_2(\bR)$ realised as irreducible   subalgebra $$\gsl_2^5(\bR)  \subset \gso(2,3)=\gso(W)$$  where $W = S^4(\bR^2)$, and as irreducible subalgebra  $$\gsl_2^4(\bR) \subset  \gsp_2(\bR) = \gsp(V)$$  where  $V = S^3 (\bR^2)$.

The remaining cases are the reductive maximal subalgebras
\begin{itemize}
	\item[(iv)] $\gso_3(\bR) \oplus \gso_2(\bR)\cong \gu_2$,
	\item[(v)] $\gso(1,2)  \oplus \gso_2(\bR)\cong   \gsl_2(\bR)^d  \oplus\gso_2(\bR)= \gu_{1,1}$.
\end{itemize}
The subalgebra $\gs_2:=\gsp_2(\bR)\cap\stab(J)$ consisting of the symplectic matrices preserving  the  complex  structure $J : p_j \mapsto q_j,\,\,  j=1,2$, is the the maximal compact subalgebra $\gu_2$ of $\gsp_2(\bR)$. We note  that  the complexification  of $\gs_2$  is  the  complex Lie    algebra
 $$ \gs_2^{\bC} \cong \ggl(V^{10})=  \ggl_2(\bC) $$
where  $V^{10}, V^{01}  \subset V^{\bC}$ are   $\pm i$ eigenspaces  of  the  complex  structure  $J$  in  $V^{\bC}$. The maximal reductive subalgebra $\gs_3:=\gu_{1,1}$ is similarly described, using a ``split'' complex structure $\widetilde J$.
\br{\rm
In case (v),  the  symplectic  space  is given by $ V  = H \otimes E = \bR^2 \otimes \bR^2$  with symplectic  form $\Omega = \mathrm{vol}_H \otimes \eta_E$,   where    $\mathrm{vol}_H$  is  the   volume  form in  $H$   and $\eta_E$ a metric in  $E$
 of  signature  $(2,0)$.  In  this  case   $ \gsl_2(\bR) \otimes \mathrm{id}=  \gsl_2(\bR)^d$  is  the  diagonal  subalgebra   of  $\gs_1$.}
\er
\br{\rm
If we were to take a metric in $E$ of signature $(1,1)$, then the associated subalgebra
would not be maximal in $\gso(2,3)$, as it stabilizes an isotropic line.}
\er
\subsection{Type  of  maximal  subalgebras   of  $\gsp_2(\bR)$}
\label{sec:mostFT}
Let $$(V, \O)=  (\bR^4, \O = -(p^*_1 \wedge q^*_1 + p^*_2 \wedge q^*_2))$$
be the $4$-dimensional symplectic vector space with symplectic basis $(p_1,p_2, q_1,q_2)$.
Using  the symplectic form, we may identify  the  symplectic Lie  algebra $\gs=\gsp(V)$  with the space
$S^2(V)$ of quadratic polynomials on $V^*$: $$uv : w \mapsto \Omega(u,w)v + \Omega(v,w)u\;,$$
for all $u,v,w\in V$. In particular $q_1p_1: p_1 \mapsto p_1, q_1 \mapsto -q_1$. More generally $\gs^{(j)} = S^{j+2}(V)$ and the full
algebra of formal symplectic vector fields is the symmetric algebra of $V$,
$$\gs^{(\infty)}  = S^+(V) = V + S^2(V) + S^3(V) + \cdots\;,$$
modulo constants.
\\

There  are two   gradings of $V$  which define two maximal  parabolic subalgebras of $\gs$. They are described by painted Dynkin diagrams, where the black node corresponds to the simple root of degree one and the white node to the simple root of degree zero.
\begin{itemize}
 \item[1.] $\bullet => \circ \qquad$  The  grading   element $d= \tfrac{1}{2}(q_1p_1 + q_2p_2)$
defines  a grading
\begin{equation*}
\begin{aligned}
V&=V^{-\tfrac{1}{2}}+V^{+\tfrac{1}{2}}\\
&= Q + P
\end{aligned}
\end{equation*}
of $V$, where $V^{-\tfrac{1}{2}}=  Q := \mathrm{span}(q_1,q_2)$, $V^{+\tfrac{1}{2}}=P := \mathrm{span}(p_1,p_2)$. The induced grading
\begin{equation*}
\begin{aligned}
\gs &=\mathfrak{s}^{-1}+ \mathfrak{s}^0 +\mathfrak{s}^1\\
&=S^2(Q) + Q\vee P + S^2(P)
\end{aligned}
\end{equation*}
of $\gs$ has corresponding parabolic  subalgebra
$\mathfrak{p}_1 = \mathfrak{s}^{\geq 0}
=  Q \vee P + S^2(P) \cong \ggl(P) + S^2(P)$. Note  that   the    associated  flag manifold $LG(V) = \mathrm{Sp}(V)/\mathrm{GL}(P)\cdot S^2 (P)$
is  the   Lagrangian  Grassmannian.
\item[2.] $\circ => \bullet  \qquad $  The grading  element $d=q_1p_1$
defines a grading
\begin{equation*}
\begin{aligned}
V&= V^{-1} + V^0 + V^1\\
&=\bR q_1 + W + \bR p_1
\end{aligned}
\end{equation*}
of $V$, where  $W=V_2 := \mathrm{span} (p_2,q_2)$. It induces the
depth  2 grading
\begin{equation*}
\begin{aligned}
\gs&=\mathfrak{s}^{-2}+ \mathfrak{s}^{-1}+ \mathfrak{s}^0+ \mathfrak{s}^1 + \mathfrak{s}^2\\
&=\RR q_1^2 + q_1 W + (S^2(W) + \RR p_1q_1) +
p_1 W + \RR p_1^2
\end{aligned}
\end{equation*}
of $\gs$ with associated parabolic subalgebra
$\mathfrak{p}_2 = \mathfrak{s}^{\geq 0} = (S^2(W) +\RR p_1q_1)+ p_1 W + \RR p_1^2$. This is a subalgebra of the  Lie  algebra of derivations of the 3-dimensional
Heisenberg Lie algebra $\mathfrak{heis}(W)= p_1 W + \RR p_1^2$.
\end{itemize}
\bt
\label{thm:maxsub}
Any maximal  subalgebra of $\gs=\gsp(V)$ is  conjugated  to one of the seven subalgebras
\begin{equation}
\begin{aligned}
\mathfrak{s}_1 &= \gsp(V_1) \oplus \gsp(V_2) \cong \gsl_2(\bR)\oplus \gsl_2(\bR)\,,\\
\gs_2 &= \gsp(V)\cap \stab(J)\cong \gu_2\,,\\
\gs_3&= \gsp(V)\cap \stab(\widetilde J)\cong \gu_{1,1}\;,\\
\gs_4&=\gsl_2(\bC)\;,\\
\gs_5&=\gsl_2^4(\bR)\;,\\
\gp_1 &=  Q\vee P + S^2(P) \cong \mathfrak{gl}_2(\bR) + S^2(\bR^2)\,,\\
\mathfrak{p}_2 &= (S^2(W)+\RR p_1q_1) + p_1 W + \bR p_1^2\cong \ggl(W)+ \mathfrak{heis}(W)\,.
\end{aligned}
\end{equation}
\et
\bc
\label{cor:finitetype}
All maximal subalgebras  of $\gs=\gsp(V)$ are  algebras of infinite  type except the Lie algebras of skew-Hermitian matrices $\gs_2$, $\gs_3$ and the singular subalgebra $\gs_5$.
\ec

This result follows from the  observation   that  their  complexifications  have  a rank one  endomorphism   and    the following  criterion, see \cite{Wil} and also \cite[Lemma 3.4]{AS}.
 \bt
\label{thm:rank1}
\hfill\par\noindent
\begin{itemize}
\item[1.] A complex linear Lie  algebra $\gh \subset \mathfrak{gl}_n(\bC)$ has infinite  type if and only if  it has  an element of  rank one.
\item[2.]
A real linear Lie  algebra $\gh \subset \mathfrak{gl}_n(\RR)$ has  infinite  type if and only if its  complexification $\gh^{\bC} \subset  \mathfrak{gl}_n(\bC)$ has infinite  type.
\end{itemize}
\et

The following criteria are for the prolongation to be trivial.
\begin{Prop}
\label{prop:fullproltrivial}
\hfill\par\noindent
\begin{itemize}
\item[1.] The first prolongation of any compact subalgebra $\gh\subset\mathfrak{so}_{n}(\bR)$ is trivial, i.e., $\gh^{(1)}=0$.
\item[2.]  Let $\gh \subset \mathfrak{gl}_n(\mathbb{K})$ be  a linear Lie  algebra over $\mathbb{K =\bC}$ or $\bR$   and  $$\gh^d_{\theta} = \{ (A,\theta A)\mid A \in \gh   \}\subset \mathfrak{gl}_n(\mathbb K) \oplus \mathfrak{gl}_n(\mathbb K)$$  the diagonal  subalgebra  of  $\mathfrak{gl}_{2n} (\mathbb K)$ twisted  by a Lie algebra automorphism $\theta:\gh\to\gh$ of $\gh$.
Then     the  first prolongation of $\gh^d_\theta$ is trivial.
\end{itemize}
\end{Prop}
In particular $\gs_2\cong\gu_2\subset \gso_4(\RR)$ satisfies $\gs_2^{(j)}=0$ for all $j\geq 1$ and likewise for $\gs_3\cong\gu_{1,1}$.
\begin{Prop}
The prolongation of $\gs_5=\gsl_2^4(\bR)$ is trivial, i.e., $\gs_5^{(j)}=0$ for all $j\geq 1$.
\end{Prop}
\pf The representation of $\gs_5$ on $S^3(\bR^2)$ is irreducible and has finite type. By a classical result of Kobayashi and Nagano, either the first prolongation is trivial or the full prolongation of $\gs_5$ is a simple Lie algebra \cite{KoNa}. In the latter case, the zero-degree component of the full prolongation must have a non-trivial center, which is not our case.
\qed
It is known that any proper subalgebra of the Lorentz algebra $\gso(1,3)$ preserves a line, either isotropic or not. This implies the following.
\begin{Prop}
\label{prop:S4}
Let $\gh$ be a proper subalgebra of $\gs_4$. Then either $\gh\subset \gs_2$ or $\gh\subset\gp_1$.
\end{Prop}

\subsection{Full prolongation  of maximal parabolics $\gp_1$, $\gp_2$ and of $\gs_1$}
In view of the results of \S \ref{sec:mostFT}, the description of the finite type subalgebras of $\gsp(V)$ reduces to the description of the finite type subalgebras of $\gs_1$, $\gp_1$ and $\gp_2$. We now determine the full prolongation of these Lie algebras. Alternative descriptions in terms of vector fields on the real plane are given in \S\ref{sec:K}.

We depart with $\mathfrak{p}_1 = Q\vee P + S^2(P)$. Its full prolongation
$$\mathfrak{p}_1^{(\infty)}=V+\mathfrak{p}_1+\mathfrak{p}_1^{(1)}+\cdots$$ is infinite-dimensional and
\begin{equation}
\begin{aligned}
\gp_1^{(k)}&=\left\{X\in\gsp(V)^{(k)}\mid\ad^{k}_V(X)\in\gp_1\right\}\\
&=Q\vee S^{k+1} (P)+S^{k+2} (P)\;,
\end{aligned}
\end{equation}
for all $k\geq 1$.
\begin{Lem}
\label{lem:fullprol1}
A subalgebra $\gh$ of $\gp_1$ has finite type if and only if $\ga=\gh\cap S^2(P)$ has finite type.
\end{Lem}
\pf
Let $\ga^{(i)}=0$ for all $i$ bigger than or equal a nonnegative integer $N$ and decompose any element $X\in\gh^{(N+1)}$ into
$X=X'+X''$, where
$$X'\in Q\vee S^{N+2}(P)\;\;,\qquad X''\in S^{N+3}(P)\;\;.$$
We have $$\ad_{P}X'=\ad_{P}X\subset \gh^{(N)}\cap S^{N+2}(P)\;,$$
that is $\ad_{P}X'\subset \ga^{(N)}=0$,
whence $X'=0$ and $X''\in \ga^{(N+1)}$. We arrive at $X=0$, so $\gh$ has finite type, proving an implication. The converse implication is trivial.
\qed

The full prolongation
$$\mathfrak{p}_2^{(\infty)}=V+\mathfrak{p}_2+\mathfrak{p}_2^{(1)}+\cdots$$
of $\mathfrak{p}_2$ is also infinite-dimensional and
\begin{equation}
\label{eq:fullprolpar2}
\begin{aligned}
\gp_2^{(k)}&=\left\{X\in\gsp(V)^{(k)}\mid\ad^{k}_V(X)\in\gp_2\right\}\\
&=\RR p_1^{k+1}q_1+ \sum_{i+j=k+2}S^{i}(W) \vee p_1^{j}\;,
\end{aligned}
\end{equation}
for all $k\geq 1$.
\begin{Lem}
\label{lem:fullprol2}
A subalgebra $\gh$ of $\gp_2$ has finite type if and only if the ideal $\widetilde\gh=\gh\cap (S^2(W)+p_1W+\RR p_1^2)$ has finite type.
\end{Lem}
\pf
Let ${\widetilde\gh}^{(i)}=0$ for all $i$ bigger than or equal a nonnegative integer $N$ and decompose any element
$X\in\gh^{(N+1)}$ into $X=X'+X''$, where
$$X'=\lambda p_1^{N+2}q_1\;\;,\qquad  X''\in\sum_{i+j=N+3}S^{i}(W) \vee p_1^{j}\;\;,$$
with $\lambda\in\RR$. Note that $\ad_{p_1}X=-\lambda p_1^{N+2}\in\widetilde\gh^{(N)}=0$, whence $X'=0$ and $X''\in\widetilde\gh^{(N+1)}=0$.
\qed
The proofs of the following two results are immediate.
\begin{Lem}
\label{prop:fullprol3}
The full prolongation of $\gs_1=\gsp(V_1)\oplus\gsp(V_2)$ is the Lie algebra direct sum of the full prolongations of $\gsp(V_1)$ and $\gsp(V_2)$.
\end{Lem}

\begin{Lem}
\label{lem:fullprol3}
 A subalgebra $\mathfrak{h}$ of $\mathfrak{s}_1$ has finite type if and only if
both $\gh\cap \gsp(V_1)$ and $\gh\cap\gsp(V_2)$ have finite type.
\end{Lem}

\section{Classification  of  finite  type  subalgebras  of  $\gsp_2(\bR)$}
\label{sec:classification}
\subsection{Finite  type  subalgebras of  \texorpdfstring{$\gs_1=\fsp(V_1)\oplus\fsp(V_2)$}{} }
\label{sec:classification3.1}
Let $V=V_1 + V_2$, $V_i=\span(p_i,q_i)$ for $i=1,2$, be  an orthogonal  decomposition of $V$. The  problem  of  description of  finite  type subalgebras of $\gs_1=\gsp(V_1)\oplus\gsp(V_2)\cong\fsl_2(\bR)\oplus\fsl_2(\bR)$ reduces to the description of the finite  type  subalgebras of $\gsl_2(\bR)$, see Lemmata \ref{prop:fullprol3} and \ref{lem:fullprol3}. We recall that the non-trivial subalgebras of $\gsl_2(\bR)$ are, up to conjugation, the diagonal subalgebra $\bR\diag(1,-1)$, $\fso_2(\bR)$, the Borel subalgebra $\gb_2$ (i.e. upper
triangular subalgebra of $\gsl_2(\bR)$) and its nilradical $\gn_2$ of strictly upper triangular matrices.

\bp  A non-zero finite type  subalgebra of $\gsl_2(\bR)$ is a Cartan subalgebra,
i.e., it is conjugated either to $\bR\mathrm{diag}(1,-1)$ or to $\fso_2(\bR)$.
\ep

\pf The only non-zero subalgebra of the complexification $\gsl_2(\bC)$ of $\gsl_2(\bR)$ with  no  rank one matrices is  the  Cartan  subalgebra. The claim  follows from this and Theorem  \ref{thm:rank1}.
\qed

We get:

\bt
\label{thm:3}
The  maximal finite  type  subalgebras  of $\gs_1\cong \fsl_2(\bR)\oplus\fsl_2(\bR)$  are (up to conjugation) the   subalgebras
\begin{itemize}
\item[1.] $\gh_1 =\fk\oplus\fk'$\;,
\item[2.] $\gh_2 =\gsl_2(\bR)^d$\;,
\item[3.] $\gh_3= \gsl_2(\bR)^d_\theta$\;,
\end{itemize}
where $\fk$, $\fk'$ are Cartan subalgebras of $\fsl_2(\bR)$, $\gh_2$ is the diagonal  subalgebra  and $\gh_3$ is the diagonal subalgebra twisted by $\theta=\Ad_{\tiny\begin{pmatrix}1& 0\\ 0 & -1\end{pmatrix}}$. The first prolongation is trivial in all cases.
\et

To prove the theorem, we first need to recall the description of subalgebras of the direct sum $\gg_1\oplus\gg_2$ of two Lie algebras. This is the content of Goursat's Lemma, see e.g. \cite{CV}.
\bl[Goursat's Lemma]
\label{goursatL}
Let $\gg_1$, $\gg_2$ be Lie algebras. There is a one-to-one correspondence between Lie subalgebras $\fh\subset \gg_1\oplus \gg_2$ and quintuples ${\mathcal Q}(\fh)=(A,A_0,B,B_0,\theta)$, with $A\subset\gg_1$, $B\subset\gg_2$ Lie subalgebras, $A_0\subset A$, $B_0\subset B$ ideals and $\theta: A/A_0\to B/B_0$ a Lie algebra isomorphism.
\el

For any subalgebra $\fh$ of $\gg_1\oplus\gg_2$ the associated quintuple ${\mathcal Q}(\fh)=(A,A_0,B,B_0,\theta)$ can be described as follows. Let $\pi_i:\gg_1\oplus\gg_2\to\gg_i$, $i=1,2$, be the natural projections and set
$$A:=\pi_1(\fh)\subset \gg_1\;,\qquad A_0:=\ker({\pi_2}_{|\fh})\cong \fh\cap\gg_1\;,$$
$$B:=\pi_2(\fh)\subset \gg_2\;,\qquad B_0:=\ker({\pi_1}_{|\fh})\cong \fh\cap\gg_2\;.$$
It is not hard to see that $A_0$ and $B_0$ can be identified with ideals in $A$ and $B$ respectively. The map $\theta:A/A_0\to B/B_0$ is defined as follows. Given $a\in A$, take any $b\in B$ such that $a+b\in \fh$ and set $\theta(a+A_0):=b+B_0$. It is easy to check that this map is well defined and it is a Lie algebra isomorphism. This gives $\fh\mapsto{\mathcal Q}(\fh)$.

Conversely, a quintuple $Q=(A,A_0,B,B_0,\theta)$ as above defines a subalgebra $\fh=\mathcal H(Q)$ by
$$\fh:=\Big\{ a+b\in A\oplus B \mid  \theta(a+A_0)=b+B_0\Big\}\;. $$
The operation ${\mathcal H}$ and ${\mathcal Q}$ are inverse to each other.
\br
{\rm We will consider different subalgebras $\fh=\mathcal H(A, A_0, B, B_0,\theta)$ for some fixed $A, A_0$ and $B, B_0$. In this case it is convenient to identify the quotients $A/A_0$ and $B/B_0$ with a fixed abstract Lie algebra $\mathfrak f$ and look at $\theta:A/A_0\to B/B_0$ as an automorphism of $\mathfrak f$.}
\er

The conjugacy classes of subalgebras are described in terms of the adjoint action of Lie groups $G_1$, $G_2$, with Lie algebras $Lie(G_1)=\gg_1$, $Lie(G_2)=\gg_2$,
on the associated quintuples.

\bl
Two subalgebras $\fh,\fh'$  of $\gg$ with corresponding quintuples $\mathcal{Q}(\fh)=(A,A_0,B,B_0,\theta)$ and
$\mathcal{Q}(\fh')=(A',A'_0,B',B'_0,\theta')$
are conjugate if and only if  there exists $(g_1,g_2)\in G_1 \times G_2$ such that
$A'=\Ad_{g_1}(A), B'=\Ad_{g_2}(B), A'_0=\Ad_{g_1}(A_0), B'_0=\Ad_{g_2}(B_0)$ and the diagram
\begin{center}
\begin{tikzpicture}[node distance=2cm, auto]
  \node (A) {$A/A_0$};
\node (B) [right of=A] {$B/B_0$};
\node (C) [below of= A] {$A'/A'_0$};
\node (D) [below of= B] {$B'/B'_0$};
---------
\draw[->] (A) to node [above] {$\theta $} (B);
\draw[->] (A) to node [swap] {$\Ad_{g_1}$} (C);
\draw[->] (B) to node [right]   {$\Ad_{g_2}$} (D);
\draw[->] (C) to node [swap] {$\theta'$} (D);
  \end{tikzpicture}
  \end{center}
  commutes.
\el

We now turn to the proof of Theorem \ref{thm:3}. We consider quintuples $Q=(A,A_0,B,B_0,\theta)$,
with $A,B$ subalgebras of $\fsl_2(\RR)$, $A_0,B_0$ ideals of $A$ and $B$, respectively,
and $\theta:A/A_0\to B/B_0$ a Lie algebra isomorphism. By Lemma \ref{lem:fullprol3}, the associated subalgebra $\fh=\mathcal H(Q)$ of $\gs_1$ has finite type if and only if $A_0$ and
$B_0$ have finite type, i.e., they are Cartan subalgebras of $\fsl_2(\bR)$ (and in that case $\gf=0$) or zero.

Table \ref{tab:G1} below describes all the possibilities up to conjugation. Subalgebras $\fh=\mathcal H(Q)$ with trivial $\gf$ are all contained in the subalgebra $\gh_1$ of Theorem \ref{thm:3}.
If $\gf$ is a Cartan subalgebra, $\theta:\gf\to\gf$ is any non-zero multiple of the identity and the twisted diagonal $\gf^d_\theta$ is again a subalgebra of $\gh_1$. If $\gf=\fn_2$, $\fb_2$ or $\fsl_2(\bR)$ then
$$\fh=\ff^d_\theta=\Big\{ a+b\in \ff\oplus \ff\mid \theta(a)=b\Big\}$$
is the  diagonal subalgebra twisted by (the restriction of) the adjoint action $\theta=\Ad_{g_1}$ of a $2\times 2$-matrix $g_1$ with determinant $\pm 1$. The corresponding subalgebras are conjugated either to subalgebras of $\gh_2$ or to subalgebras of $\gh_3$ as in Theorem \ref{thm:3}.
\vskip0.5cm\par
{
\begin{table}[H]
\footnotesize
\begin{centering}
\makebox[\textwidth]{%
\begin{tabular}{|c|c|c|}
\hline
$(A,A_0)$ & $(B,B_0)$ & $\gf$\\
\hline
\hline
$\begin{gathered}\\ A=A_0=\{0\},\;\bR\diag(1,-1)\;\text{or}\;\fso_2(\bR)\\ \\ \end{gathered}$ & $\begin{gathered}\\ B=B_0=\{0\},\;\bR\diag(1,-1)\;\text{or}\;\fso_2(\bR)\\ \\ \end{gathered}$ & $\{0\}$ \\
\hline
$\begin{gathered}\\ A=\bR\diag(1,-1),\; A_0=\{0\}\\ \\ \end{gathered}$ & $\begin{gathered}\\ B=\bR\diag(1,-1),\; B_0=\{0\}\\ \\ \end{gathered}$ & $\bR\diag(1,-1)$\\
\hline
$\begin{gathered}\\ A=\fso_2(\bR),\;A_0=\{0\}\\ \\ \end{gathered}$ & $\begin{gathered}\\ B=\fso_2(\bR),\;B_0=\{0\}\\ \\ \end{gathered}$ & $\fso_2(\bR)$\\
\hline
$\begin{gathered}\\ A=\fn_2,\;A_0=\{0\}\\ \\ \end{gathered}$ & $\begin{gathered}\\ B=\fn_2,\;B_0=\{0\}\\ \\ \end{gathered}$ & $\mathfrak n_2$\\
\hline
$\begin{gathered}\\ A=\fb_2,\;A_0=\{0\}\\ \\ \end{gathered}$ & $\begin{gathered}\\ B=\fb_2,\;B_0=\{0\}\\ \\ \end{gathered}$ & $\mathfrak b_2$\\
\hline
$\begin{gathered}\\ A=\fsl_2(\bR),\;A_0=\{0\}\\ \\ \end{gathered}$ & $\begin{gathered}\\ B=\fsl_2(\bR),\;B_0=\{0\}\\ \\ \end{gathered}$ & $\fsl_2(\bR)$\\
\hline
\end{tabular}}
\end{centering}
\caption[]{\label{tab:G1} Goursat quintuples for finite type subalgebras of $\fsl_2(\bR)\oplus\fsl_2(\bR)\subset \fsp_2(\bR)$.} \vskip14pt
\end{table}}
\vskip-0.4cm\par
\subsection{Finite  type  subalgebras of \texorpdfstring{$\mathfrak{p}_1 = QP + S^2(P)$}{}}
The description of the finite  type subalgebras of $\gp_1$ reduces to the description of finite  type  subspaces of $S^2(P)$, see Lemma \ref{lem:fullprol1}. Note that $ S^2 (P)  = \mathrm{span}(p_1^2, p_1p_2, p_2^2) \subset \gsp(V) $  has a rank one   endomorphism
and  it is therefore of infinite  type. On the other hand $QP\cong\ggl(P)$ acts diagonally on $$V=P+Q\cong P+P^*$$ and its first prolongation is trivial, see Proposition \ref{prop:fullproltrivial}.
 \bl
\label{lemma:lineinfinite}
The $1$-dimensional subspace of $S^2(P)$ generated by $X=x_1 p_1^2 + x_2 p_1p_2 + x_3 p_2^2$ has  infinite  type if  and only if $x_2^2 =4x_1 x_3$.
   \el
\pf After complexification, we  calculate
\begin{equation*}
\begin{aligned}
X &: q_1 \mapsto -2x_1 p_1 -x_2 p_2 \mapsto 0\;,\\
X &: q_2 \mapsto -2x_3 p_2 - x_2 p_1 \mapsto 0\;,
\end{aligned}
\end{equation*}
and note that  $X$  has rank one if  and only if   $X(q_1)$  and  $X(q_2)$ are proportional.
 \qed
Lemma \ref{lemma:lineinfinite} admits a more suggestive interpretation. Under the identification of $\mathrm{SL}(P)$ with the (connected) spin group $\mathrm{Spin}^\circ(1,2)$ in Lorentzian signature, the action on $S^2(P)$  corresponds to the vectorial representation of $\mathrm{Spin}^\circ(1,2)$ on the Lorentzian vector space $\bR^{1,2}$ and the elements of infinite type to the lightlike vectors.
We shall fix an orthonormal basis $(\Be_0,\Be_1,\Be_2)$ of $S^2(P)$,
\begin{equation}
\label{eq:basislor}
\Be_0=\tfrac{1}{2}(p_1^2+p_2^2)\;,\qquad\Be_1=\tfrac{1}{2}(p_1^2-p_2^2)\;,\qquad \Be_2=p_1p_2\;,
\end{equation}
so that {\it finite type} $1$-dimensional subspaces of $S^2(P)$ are conjugated either to $\bR\Be_0$ (positive norm) or $\bR\Be_2$  (negative norm). One can directly check that the first prolongation is trivial in both cases, i.e., $(\bR\Be_0)^{(1)}=(\bR\Be_2)^{(1)}=0$.
\bl
\label{lem:subspaceinf}
Any $2$-dimensional subspace of  $S^2 (P)$ has infinite type.
 \el
\pf
After complexification, it includes a lighlike vector.
\qed
A direct consequence of
Lemma \ref{lem:subspaceinf} and the proof of Lemma \ref{lem:fullprol1} is the following a priori estimate of the type of $\gh$:
\bc
Let $\gh\subset \fp_1$ be a finite type subalgebra. Then $\dim(\gh\cap S^2(P))\leq 1$ and $\gh^{(k)}=0$ for all $k\geq 2$.
\ec

To proceed further, we shall use the classification of the subalgebras of $\fp_1$. Note first that $\fp_1\cong\mathfrak{sim}(1,2)=\mathfrak{co}(1,2)+\bR^{1,2}$, $QP\cong\mathfrak{co}(1,2)$, $S^2(P)\cong \bR^{1,2}$, is isomorphic to the similitude algebra of the Lorentzian vector space $\bR^{1,2}$.
The description of the conjugacy classes of subalgebras of the similitude algebra  is rather involved and it can be found in \cite[$\S 4$, Tables $2$-$4$]{PWSZ}.
We recall here only the facts that we need and refer directly to \cite{PWSZ} for more details (see also the overview \cite{Win}).
\vskip0.4cm\par

Let  $\mathfrak{h} \subset \mathfrak{sim}(1,2)$  be a subalgebra, $\overline{\mathfrak{h}} $  its  projection to  $\mathfrak{co}(1,2)\cong\mathfrak{sim}(1,2)/\bR^{1,2}$
 and    $\mathfrak{a} = \mathfrak{h} \cap \bR^{1,2}$. The subalgebra $\mathfrak h_s=\overline{\mathfrak h}+\mathfrak a$ is called the {\it splitting subalgebra associated to $\fh$} and it is an invariant of $\mathfrak h$ up to conjugation.  All subalgebras conjugated to a subalgebra of  the   form  $\mathfrak{h}= \overline{\mathfrak{h}}+\mathfrak{a}$ are called {\it splitting}; otherwise they are {\it nonsplitting}.

Nonsplitting subalgebras are all obtained as deformations of splitting subalgebras $\mathfrak{h}=\overline{\mathfrak{h}}+\mathfrak{a}$. More precisely $\overline{\mathfrak h}$ is twisted
by cocycles $Z^{1}(\overline{\mathfrak h},\mathfrak m)$ of the Chevalley-Eilenberg complex of $\overline{\mathfrak h}$ with values in $\mathfrak m=\bR^{1,2}/\mathfrak a$, and coboundaries give subalgebras conjugated under the translation subgroup $\exp(\bR^{1,2})$. In other words nonsplitting subalgebras $\fh$ with underlying associated splitting subalgebra $\fh_s=\overline{\mathfrak h}+\ga$ are described by elements of  the first Chevalley-Eilenberg cohomology group $H^{1}(\overline{\mathfrak h},\mathfrak m)$.
We note that $$\ga=\fh\cap\bR^{1,2}=\mathfrak h_s\cap\bR^{1,2}$$ so that $\gh$ {\it has finite type if and only if} $\gh_s$ {\it has finite type}, cf. Lemma \ref{lem:fullprol1}.

We make contact with the notation of \cite{PWSZ} and introduce the dilation, boost and rotation matrices
\begin{equation}
\label{eq:basislormat}
\begin{aligned}
F&=-\frac{1}{2}(p_1q_1+p_2q_2)\;,\qquad
K_1=\frac{1}{2}(p_1q_1-p_2q_2)\;,\\
K_2&=\frac{1}{2}(p_1q_2+p_2q_1)\;,\qquad\;\;\;
L_3=-\frac{1}{2}(p_1q_2-p_2q_1)\;,
\end{aligned}
\end{equation}
which together with \eqref{eq:basislor} identify $\fp_1$ with $\mathfrak{sim}(1,2)$. The brackets are as in \cite[p. 951]{PWSZ}, in particular $[F,\Be_i]=-\Be_i$ for $i=0,1,2$. The following result follows from Lemmata \ref{lem:fullprol1}, \ref{lemma:lineinfinite} and \ref{lem:subspaceinf},
and the classification of the subalgebras of $\mathfrak{sim}(1,2)$.
\bp
The finite type subalgebras $\gh$ of $\gp_1$ are, up to conjugation, either contained in $QP\cong\mathfrak{co}(1,2)$ or they are one of the subalgebras listed in Table \ref{tab:G2}:
 {
\begin{table}[H]
\small
\begin{centering}
\makebox[\textwidth]{%
\begin{tabular}{|c|c|c|c|}
\hline
Name & Abstract Lie algebra & Generators & Splitting Subalgebra\\
\hline
\hline
$F_{6,5}$ & $\bR$ & $\Be_2$ &\\
\hline
$F_{6,6}$ & $\bR$ & $\Be_0$ &\\
\hline
$F_{3,5}$ & $\fso(1,1)\oplus\bR$ & $K_1,\Be_2$ &\\
\hline
$F_{5,3}$ & $\fso_2(\bR)\oplus\bR$ & $L_3,\Be_0$ &\\
\hline
$DF_{6,5}$ & $\bR+\bR$ & $F, \Be_2$ &\\
\hline
$DF_{6,6}$ & $\bR+\bR$ & $F, \Be_0$ &\\
\hline
$DF_{3,5}$ & $\mathfrak{co}(1,1)+\bR$ & $F, K_1,\Be_2$ &\\
\hline
$DF_{5,3}$ & $\mathfrak{co}(2)+\bR$ & $F, L_3,\Be_0$ &\\
\hline
$\tilde{F}_{3,9}^{\phantom{A^A}}$ & $\bR$ & $K_1+a\Be_2$\;\;\text{where}\;\;$a\neq 0$ & $\bR K_1$\\
\hline
$\tilde F_{4,7}^{\phantom{A^A}}$ & $\bR$ & $K_2+L_3+\epsilon(\Be_0+\Be_1)$\;\;\text{where}\;\;$\epsilon=\pm 1$ & $\bR(K_2+L_3)$\\
\hline
$\tilde F_{5,6}^{\phantom{A^A}}$ & $\bR$ & $L_3+a\Be_0$\;\;\text{where}\;\;$a\neq 0$ & $\bR L_3$\\
\hline
$D_{4,12}$ & $\bR+\bR$ & $F-K_1+\epsilon(\Be_0-\Be_1), K_2+L_3$\;\;\text{for}\;\;$\epsilon=\pm 1$ & $D_{4,11}$\;\;\text{with}\;\;$a=-1$\\
\hline
$D_{4,13}$ & $\bR+\bR$ & $F+(1/2)K_1, K_2+L_3+\epsilon(\Be_0+\Be_1)$\;\;\text{for}\;\;$\epsilon=\pm 1$ & \,$D_{4,11}$\;\;\text{with}\;\;$a=1/2$\\
\hline
$D_{6,13}$ & $\bR+\bR$ & $F+aK_1, \Be_2$\;\;\text{for}\;\;$a>0$ &\\
\hline
$D_{6,14}$ & $\bR+\bR$ & $F+K_1+\epsilon(\Be_0+\Be_1), \Be_2$\;\;\text{for}\;\;$\epsilon=\pm 1$ & $D_{6,13}$\;\;\text{with}\;\;$a=+1$\\
\hline
$D_{6,15}$ & $\bR+\bR$ & $F+aL_3, \Be_0$\;\;\text{for}\;\;$a\neq0$ &\\
\hline
$D_{6,22}$ & $\bR$ & $F+K_1+\epsilon(\Be_0+\Be_1)$\;\;\text{for}\;\;$\epsilon=\pm 1$ & $D_{6,21}$\;\;\text{with}\;\;$a=+1$\\
\hline
\end{tabular}}
\end{centering}
\caption[]{\label{tab:G2} Finite type subalgebras of $\gp_1\subset\fsp_2(\bR)$ not contained in $QP$.} \vskip14pt
\end{table}}
\vskip-0.7cm\par\noindent
The associated splitting subalgebras have been illustrated for all the nonsplitting subalgebras and we used the notation $D_{4,11}:=\span(F+aK_1,K_2+L_3)$ and $D_{6,21}:=\span(F+aK_1)$.
\ep

The  maximal finite type subalgebras of $\fp_1$ are $\mathfrak{co}(1,2)$, $DF_{3,5}$, $DF_{5,3}$
and the nonsplitting subalgebras $D_{4,12}$, $D_{4,13}$ and $D_{6,14}$. First note that the  commutative radical $\bR^{1,2}$ of $\fp_1$ is an irreducible $\mathfrak{co}(1,2)$-module, whence $\mathfrak{co}(1,2)$ is a maximal subalgebra. The subalgebras $DF_{3,5}$ and $DF_{5,3}$ are $3$-dimensional and cannot be
(conjugated in $\gp_1$ to) subalgebras of $\mathfrak{co}(1,2)$, as they intersect $\bR^{1,2}$ nontrivially. Finally, the subalgebras $D_{4,12}$, $D_{4,13}$ and $D_{6,14}$ are all $2$-dimensional and clearly not subalgebras of $\mathfrak{co}(1,2)$. A straightforward computation shows that they cannot be included in any $3$-dimensional finite type subalgebra of $\gp_1$ either, in particular they are not conjugated to subalgebras of $DF_{3,5}$ and $DF_{5,3}$.
\\

We have proved most of the following.
\bt
\label{thm:4}
The  maximal finite  type  subalgebras of $\gp_1\cong \mathfrak{co}(1,2)+\bR^{1,2}$  are, up to conjugation, given by
\begin{itemize}
\item[1.] $\fh=\mathfrak{co}(1,2)$,
\item[2.] the $3$-dimensional solvable splitting subalgebras $\fh=\overline{\mathfrak h}+\mathfrak a$, where
$\mathfrak a=\bR\Be_0$ or $\bR\Be_2$ and $\overline{\mathfrak h}=\mathrm{Nor}_{\mathfrak{co}(1,2)}(\mathfrak a)$ is the normalizer of $\mathfrak a$ in $\mathfrak{co}(1,2)$,
\item[3.] the $2$-dimensional nonsplitting subalgebras $D_{4,12}$, $D_{4,13}$ and $D_{6,14}$ illustrated in Table \ref{tab:G2}.
\end{itemize}
The first prolongation is trivial in all cases.
\et
\pf
It remains to show the last claim. If $\gh=\mathfrak{co}(1,2)$ we already know that $\gh^{(1)}=0$ and the other
cases are settled by straightforward computations, which we omit.
\qed

\subsection{Finite  type  subalgebras of \texorpdfstring{$\gp_2= \ggl(W)+ \mathfrak{heis}(W)$}{} }
By Lemma \ref{lem:fullprol2}, the description of the finite  type subalgebras of
$\fp_2 =\ggl(W)+ \mathfrak{heis}(W)$
reduces to the description of the finite  type  subalgebras of its maximal ideal $$\widetilde\fp_2=\fsl(W)+ \mathfrak{heis}(W)\;.$$
We first note that the  action of the basic  endomorphisms   of  $\mathfrak{heis}(W) = \mathrm{span} (p_1p_2, p_1q_2,p_1^2)$ on the basic elements of $V$ is given by
\begin{equation*}
\begin{aligned}
p_1p_2&: q_1 \mapsto - p_2\;,\,  q_2 \mapsto -p_1\;,\\
p_1q_2&: q_1 \mapsto - q_2\;,\,  p_2 \mapsto p_1\;,\\
p_1^2&:  q_1 \mapsto - 2p_1\;,
\end{aligned}
\end{equation*}
and trivial otherwise. This shows the following.
\bl

Any finite type subalgebra of $\mathfrak{heis}(W)$ is conjugated to one of the form $\bR p_1w$ for some $w\in W$. We  may  assume  $w=0$ (trivial subalgebra) or  $w = p_2$.
\el
Assume  $\widetilde\gh $ is  a  finite  type  subalgebra of $\widetilde\gp_2$ with nontrivial intersection $\ga=\widetilde\gh\cap\mathfrak{heis}(W)$  with  $\mathfrak{heis}(W)$.
 Then  $\ga = \bR p_1 p_2$ is  an ideal of $\widetilde\gh$, that is
$$
\widetilde \gh\subset \mathrm{Nor}_{\widetilde\gp_2}(\ga)=\mathrm{span}(p_2^2,p_2q_2)+\bR p_1p_2+\bR p_1^2\;,
$$
where $\mathrm{Nor}_{\widetilde\gp_2}(\ga)$ is the normalizer of $\ga$ in $\widetilde\gp_2$. We note that
the normalizer is the semidirect sum of the Borel subalgebra $\gb_2$ of $\fsl(W)$ and a two-dimensional {\it abelian} ideal. Hence, we may regard $\widetilde\fh$ as a (possibly trivial) deformation of the corresponding splitting subalgebra $\widetilde\gh_s=\overline{\mathfrak h}+\mathfrak a$, where $\overline{\mathfrak h}$ is the  projection of $\widetilde\fh$ to $\gb_2$ (see \cite{Win}).
 The list of such splitting subalgebras consists of
\be
\label{eq:sswitha}
\begin{aligned}
\gb_2&+\bR p_1p_2\;,\qquad\qquad \gn_2+\bR p_1p_2\;,\\
\bR\diag(1,-1)&+\bR p_1p_2\quad\,\text{and}\qquad\qquad \bR p_1p_2\;,
\end{aligned}
\ee
where $\gn_2=\bR p_2^2$ is the nilradical of $\gb_2$ and $\bR\diag(1,-1)=\bR p_2q_2$ the diagonal subalgebra.
\\

The proof of the following lemma is straightforward.
\bl
The Chevalley-Eilenberg cohomology group $H^{1}(\overline{\mathfrak h},\mathfrak m)$, $\mathfrak m\cong \bR p_1^2$, is trivial when $\overline \gh$ is trivial and it coincides with the $1$-dimensional space $\mathrm{Hom}(\overline\gh,\gm)$ when $\overline\gh$ is $\gn_2$ or $\bR\diag(1,-1)$. The group $H^{1}(\gb_2,\mathfrak m)$ is $1$-dimensional too and it is generated by the cocycle $$c:\gb_2\to\mathfrak m$$ given by
$c(p_2^2)=0$ and $c(p_2q_2)=p_1^2$.
\el

It follows that the finite type subalgebras of $\widetilde\gp_2$ with nontrivial intersection with
$\mathfrak{heis}(W)$ are among the following subalgebras:
\begin{itemize}
	\item[(i)] $\mathrm{span}(p_2^2,p_2q_2+\epsilon p_1^2,p_1p_2)$,
	\item[(ii)] $\mathrm{span}(p_2^2+\epsilon p_1^2,p_1p_2)$,
	\item[(iii)] $\mathrm{span}(p_2q_2+\epsilon p_1^2,p_1p_2)$,
	\item[(iv)] $\mathrm{span}(p_1p_2)$,
\end{itemize}
where $\epsilon\in\bR$. Using conjugation by the grading element of $\gp_2$, we may also arrange for $\epsilon=0, \pm 1$.

The  action  of the basic  endomorphisms of  $\fsl(W)=\mathrm{span}(p_2^2,p_2q_2,q_2^2)$ on $W$ is
\begin{equation*}
\begin{aligned}
p_2^2=\begin{pmatrix}
0 & -2 \\
0 & 0
\end{pmatrix}\;,\quad
p_2q_2=\begin{pmatrix}
1 & 0 \\
0 & -1
\end{pmatrix}\;,\quad
q^2_2 =\begin{pmatrix}
0 & 0 \\
2 & 0
\end{pmatrix}\;,
\end{aligned}
\end{equation*}
and a finite  type subalgebra of $\fsl(W)$ consists only of semisimple  elements. This shows that the subalgebra (i) has infinite type and, in a similar way, one sees that the complexification of (ii) has an element of rank one. Using equation \eqref{eq:fullprolpar2} with $k=1$, one finally arrives at the following.
\bp
\label{prop:partial1}
The finite type subalgebras of $\widetilde\gp_2$ with nontrivial intersection with $\mathfrak{heis}(W)$ are, up to conjugation, the subalgebra $(iii)$ with $\epsilon=0$ or $\epsilon=\pm 1$ and the subalgebra $(iv)$. The first prolongation is trivial in all cases.
\ep

Assume now  that $\widetilde\gh$ has no intersection  with $\mathfrak{heis}(W)$. Let $\pi:\widetilde\fp_2\to\widetilde\fp_2/\mathfrak{heis}(W)\cong \gsl(W)$ be the natural projection and set $\overline{\mathfrak h}=\pi(\widetilde\gh)$. Clearly $\widetilde\gh\cong \overline\gh$ is isomorphic to a subalgebra of $\gsl(W)$ and we may assume $\overline\gh$ to be the whole $\gsl(W)$, $\bR\diag(1,-1)$, $\fso_2(\bR)$, $\gb_2$ or $\gn_2$. However it is {\it not} true, in general, that $\widetilde\gh$ has finite type if and only if $\overline\gh$ has finite type.

Write
$$
\widetilde\gh=\left\{(X,\varphi(X))\mid X\in\overline{\gh}\right\}\;,
$$
where $\varphi:\overline\gh\to\mathfrak{heis}(W)$. Since the Heisenberg algebra is not abelian, the map $\varphi$ has to satisfy a non-linear version of the Chevalley-Eilenberg cocycle condition in order for $\widetilde\gh$ to be a Lie algebra:
\begin{equation}
\label{eq:nonlincc}
\varphi[X,Y]=[\varphi(X),Y]+[X,\varphi(Y)]+[\varphi(X),\varphi(Y)]\;,
\end{equation}
for all $X,Y\in\overline\gh$. It is convenient to decompose $\varphi=c+\psi$ into components
$$
c:\overline\gh\to p_1W\;,\qquad \psi:\overline\gh\to \bR p_1^2\;,
$$
and rewrite \eqref{eq:nonlincc} as
\begin{equation}
\label{eq:ccdecomp}
\begin{aligned}
c[X,Y]&=[c(X),Y]+[X,c(Y)]\;,\\
\psi[X,Y]&=[c(X),c(Y)]\;.
\end{aligned}
\end{equation}
We note that $c$ is a cocycle in the Chevalley-Eilenberg complex of $\overline{\mathfrak h}$ with values in $p_1W$ and that the coboundaries correspond to subalgebras conjugated under the analytic subgroup $\exp(p_1W)$. For any given cocycle, one easily reconstructs the full map $\varphi$ using \eqref{eq:ccdecomp}.
\bl
\label{lem:CE0}
The Chevalley-Eilenberg cohomology group $H^{1}(\overline{\mathfrak h},p_1W)$ is zero,
except for $\overline\gh=\gn_2$, in which case it is $1$-dimensional and generated by the cocyle
$c(p_2^2)=p_1q_2$.
\el

It follows that the finite type subalgebras of $\widetilde\gp_2$ with trivial intersection with
$\mathfrak{heis}(W)$ are among the following subalgebras:
\begin{itemize}
	\item[(v)] $\gsl(W)$,
	\item[(vi)] $\mathrm{span}(p_2^2,p_2q_2+\epsilon p_1^2)$,
	\item[(vii)] $\mathrm{span}(p_2q_2+\epsilon p_1^2)$,
	\item[(viii)] $\mathrm{span}(p_2^2+q_2^2+\epsilon p_1^2)$,
		\item[(ix)] $\mathrm{span}(p_2^2+\epsilon p_1^2)$,
	\item[(x)] $\mathrm{span}(p_2^2+\epsilon p_1q_2)$,
\end{itemize}
where $\epsilon=0$ or $\epsilon=\pm 1$. The last two conjugacy classes correspond to the case $\overline \gh=\gn_2$ with the corresponding cocycle either vanishing or, respectively, as in Lemma \ref{lem:CE0}. The proof of the following result is similar to that of Proposition \ref{prop:partial1} and we omit it.
\bp
\label{prop:partial2}
The finite type subalgebras of $\widetilde\gp_2$ with trivial intersection with $\mathfrak{heis}(W)$ are, up to conjugation, the subalgebras $(vii)$-$(viii)$ with $\epsilon=0,\pm 1$ and the subalgebras $(ix)$-$(x)$ with $\epsilon=\pm 1$. The first prolongation is trivial in all cases.
\ep

Let us now sum up the results of Propositions \ref{prop:partial1} and \ref{prop:partial2}.
Let $\gh$ be a finite type subalgebra of $\gp_2$ and $\widetilde\gh=\gh\cap\widetilde\gp_2$ its ideal, which is a subalgebra of $\widetilde\gp_2$ as we just determined.
We either have $\gh=\widetilde\gh$ or $\widetilde\gh$ is a codimension one ideal of $\gh$.
The latter case happens if $\widetilde\gh$ is a graded subalgebra of $\widetilde\gp_2$ and $\gh=\widetilde\gh+\bR p_1q_1$, and in few more cases:
\begin{equation*}
\begin{aligned}
\mathrm{span}&(p_1p_2,p_1q_1+\lambda p_2q_2)\;,\;\qquad\quad\mathrm{span}(p_1p_2,p_1q_1+\epsilon p_2^2)\;,\\
\mathrm{span}&(p_2^2+\epsilon p_1^2,p_1q_1+p_2q_2)\;,\qquad\mathrm{span}(p_2^2+\epsilon p_1q_2,3p_1q_1+p_2q_2)\;,
\end{aligned}
\end{equation*}
where $\epsilon=\pm 1$ and $\lambda\neq 0$.
\bt
\label{thm:5}
The  maximal finite  type  subalgebras  of $\gp_2$  are, up to conjugation, given by the graded  subalgebras
\begin{itemize}
\item[1.] $\mathrm{span}(p_2q_2,p_1q_1,p_1p_2)$,
\item[2.] $\mathrm{span}(p_2^2+q_2^2,p_1q_1)$,
\end{itemize}
and the subalgebras
\begin{itemize}
\item[3.] $\mathrm{span}(p_2q_2+\epsilon p_1^2,p_1p_2)$,
\item[4.] $\mathrm{span}(p_2^2+q_2^2+\epsilon p_1^2)$,
\item[5.] $\mathrm{span}(p_2^2+\epsilon p_1^2,p_1q_1+p_2q_2)$,
\item[6.] $\mathrm{span}(p_2^2+\epsilon p_1q_2,3p_1q_1+p_2q_2)$,
\end{itemize}
where $\epsilon=\pm 1$.
The first prolongation is trivial in all cases.
\et
\pf
The first claim is now immediate and the last follows from a computation.
\qed
\br{\rm
The subalgebras $(1)$, $(3)$ and $(6)$ of Theorem \ref{thm:5} appeared already in Theorem \ref{thm:4} and the subalgebra $(2)$ in Theorem \ref{thm:3}. The subalgebra $(5)$ is included in $DF_{5,3}$ if $\epsilon=+1$ and in $DF_{3,5}$ if $\epsilon=-1$. (In the second case, it is enough to notice that both vectors $\Be_1,\Be_2$ in equation \eqref{eq:basislor} have negative norm.)}
\er
\subsection{Proof of Theorem \ref{thm:1}}
\label{sec:finalthm1}
In  this  section we shall complete the proof of Theorem \ref{thm:1}. Let us first summarize the discussion of \S\ref{sec:maximal} and Theorems \ref{thm:3}--\ref{thm:5} into the following auxiliary result. 
\bp
\label{thm:1preliminary}
A finite type subalgebra $\gh$ of $\gsp_2(\bR)$ is, up to conjugation, included in one of the subalgebras of the following list:
\begin{itemize}
\item[1.] the unitary algebra $\gu_2$;
\item[2.] the pseudo-unitary algebra $\gu_{1,1}$;
\item[3.] the subalgebra $Q\vee P\cong\ggl(P)$, where $P$ and $Q$ are  complementary  Lagrangian subspaces;
\item[4.] the irreducible subalgebra $\fsl_2^4(\bR)$ acting on $ V =  S^3(\bR^2)$;
\item[5.] the solvable   $3$-dimensional subalgebras $DF_{5,3}$ and $DF_{3,5}$;
\item[6.] the diagonal subalgebra $\fsl_2(\bR)^d$ and the twisted diagonal subalgebra $\fsl_2(\bR)^d_\theta$;
\item[7.] the direct sum $\fk\oplus\fk'$ of two Cartan subalgebras of $\fsl_2(\bR)$ (i.e., $\bR\mathrm{diag}(1,-1)$ or $\fso_2(\bR)$);
\item[8.] the solvable $2$-dimensional subalgebras $D_{4,12}$, $D_{4,13}$ and $D_{6,14}$;
\item[9.] the  $1$-dimensional subalgebra $\mathrm{span}(p_2^2+q_2^2+\epsilon p_1^2)$, where $\epsilon=\pm 1$.
\end{itemize}
\ep
We claim that all   these  subalgebras  are  maximal  finite   type   subalgebras  with  the exception  of
$$DF_{3,5},\;\; DF_{5,3},\;\;\fsl_2(\bR)^d,\;\; \fsl_2(\bR)^d_{\theta},\;\;2\bR\mathrm{diag}(1,-1),\;\;2\fso_2(\bR),\;\;  D_{4,13}.$$
Now, the $4$-dimensional subalgebras are clearly maximal finite type subalgebras of $\fsp_2(\bR)$. The complexification of the singular subalgebra $\fsl_2^4(\bR)$ is a maximal subalgebra of $\fsp_2(\bC)$ by a classical result of Dynkin \cite[Theorem 3.3]{GOVIII}, hence $\fsl_2^4(\bR)$ is a maximal subalgebra of $\fsp_2(\bR)$.

We have seen that the $3$-dimensional solvable subalgebras $DF_{5,3}$ and $DF_{3,5}$ are not
conjugated in $\gp_1$ to subalgebras of $\mathfrak{gl}(P)$. However, a simple computation shows that the ring of invariant endomorphisms of $V$ is two-dimensional, generated in the first case (resp. the second case) by a split complex structure compatible with $\Omega$ (resp. a paracomplex structure). Hence, a posteriori, $DF_{5,3}\subset\mathfrak{u}_{1,1}$ and $DF_{3,5}\subset\mathfrak{gl}(P)$, up to conjugation in the full symplectic group. 

The diagonal subalgebra $\fsl_2(\bR)^d$ is also not maximal, as $\fsl_2(\bR)^d=\mathfrak{su}_{1,1}\subset \gu_{1,1}$. Similarly, the twisted diagonal subalgebra $\fsl_2(\bR)^d_{\theta}$ preserves a paracomplex structure compatible with $\Omega$, so that $\fsl_2(\bR)^d_{\theta}=\gsl(P)\subset\ggl(P)$.

The direct sum $2\bR\mathrm{diag}(1,-1)$ of two noncompact Cartan subalgebras of $\fsl_2(\bR)$
coincides with the normalizer $\overline{\mathfrak h}=\operatorname{Nor}_{\ggl(P)}(\ga)$ of
$\mathfrak a=\bR(p_1p_2)$, hence it is contained in $DF_{3,5}\subset\mathfrak{gl}(P)$. The direct sum $2\fso_2(\bR)$ of two compact Cartan subalgebras of 
$\fsl_2(\bR)$ clearly is, up to conjugation, a subalgebra of the maximal compact subalgebra $\mathfrak u_2$. 

On the contrary, the direct sum $\bR\mathrm{diag}(1,-1)\oplus\fso_2(\bR)$ is a maximal finite type subalgebra. The ring of invariant endomorphisms of $V$ is $4$-dimensional and it is not difficult to see that there is no invariant complex or paracomplex structure compatible with $\Omega$. Hence, this Cartan subalgebra is not (conjugated to) a subalgebra of a $4$-dimensional finite type subalgebra. Furthermore, it is not contained in a $3$-dimensional finite type subalgebra either, since $\gsl_2(\bR)$ does not have $2$-dimensional subalgebras consisting of semisimple elements.

So far we obtained the following list of maximal finite type subalgebras 
\begin{equation}
\label{eq:list}
\mathfrak{u}_2,\;\;\mathfrak{u}_{1,1},\;\;\mathfrak{gl}(P),\;\;\mathfrak{sl}_{2}^4(\bR)\;\;\text{and}\;\;\bR\mathrm{diag}(1,-1)\oplus\gso_2(\bR).
\end{equation} 
We note that all such subalgebras are algebraic in the sense that, after complexification, the semisimple and nilpotent components of any element of the subalgebra still lie in the subalgebra. It remains to consider the subalgebras $8$-$9$ of Proposition \ref{thm:1preliminary}.

We depart with the subalgebras $8$ which are not conjugated to a subalgebra  from   this list. 
The algebraic  closure  of  subalgebras $D_{4,12}=\mathrm{span}(p_1 q_1+\epsilon p_2^2, p_2 q_1)$ and $D_{6,14}=\mathrm{span}(p_2q_2+\epsilon p_1^2, p_1p_2)$ are given by
\begin{equation*}
\begin{aligned}
\overline{D}_{4,12}&= \mathrm{span}(p_1 q_1 , p_2 q_1, p_2^2)\;,\\
 \overline{D}_{6,14}&=\mathrm{span}(p_2q_2 , p_1p_2, p_1^2)\;,
\end{aligned}
\end{equation*}
and they have infinite type, since they include a rank one element. Therefore $D_{4,12}$ and $D_{6,14}$ are maximal finite type subalgebras. Actually, it is not difficult to see these two subalgebras do coincide, up to conjugation in the full symplectic group. 
 
	The derived ideal of the algebra  $ D_{4,13} = \mathrm{span}(p_1q_1 + 3p_2q_2, p_2q_1 + \epsilon p_1^2) $  consists of   a nilpotent  endomorphism with Jordan   block  of  order four. Such  endomorphism   does not  exist  in  the derived ideal of any of the algebras from the list \eqref{eq:list}, with the exception of $\gsl_2^4(\bR)$. Indeed, it is easy to see that $D_{4,13}$ is exactly the Borel subalgebra $\gb_2$ of $\gsl_2^4(\bR)$.
	
Finally, the algebraic closure of the $1$-dimensional  subalgebra $9$ is
$\mathrm{span}(p_2^2+ q_2^2, p_1^2)$, which has infinite type, it is abelian and has a semisimple element generating a compact one-parameter subgroup.
One  can  easily  check  that the algebras from list \eqref{eq:list} and $\overline D_{4,12}$ do not  contain   subalgebras of this kind.
The proof of Theorem \ref{thm:1} is completed.
\subsection{Proof of Theorem \ref{thm:dichotomy}}
By the discussion above Theorem \ref{thm:dichotomy}, we already know that if $\gh$ has finite type or $\gh=\mathfrak{sp}(V)$ then $gr(\gg)=V+\gh$. Otherwise, $\gh$ is contained in a maximal infinite type subalgebra of $\mathfrak{sp}(V)$, that is a maximal parabolic subalgebra or one of the subalgebras $\fs_1$ and $\gs_4$. 

If $\fh=\fs_1$ or $\fs_4$ then, after complexification, the first prolongation $\gh^{(1)}=S^3(\mathbb C^2)\oplus S^3(\mathbb C^2)$ and by an argument similar to the case $\fh=\mathfrak{sp}(V)$ we have that $gr(\mathfrak g)$ is either infinite-dimensional, which is not possible, or $gr(\mathfrak g)=V+\fh$. Furthermore, if $\gh$ is a proper subalgebra of $\gs_4$ then either $gr(\mathfrak g)=V+\fh$ or $\fh$ is a subalgebra of $\gp_1$ by Proposition \ref{prop:S4}.

It remains to study the case where $\fh$ is a proper subalgebra of $\fs_1$ which is not of finite type. According to \S \ref{sec:classification3.1} we have
$$\fh=\Big\{ a+b\in A\oplus B \mid  \theta(a+A_0)=b+B_0\Big\}\;,$$
for some subalgebras $A,B\subset \fsl_2(\bR)$ and ideals $A\subset A_0$, $B\subset B_0$ such that $A/A_0\cong B/B_0$. First of all we can assume that
$A$ and $B$ are equal to $\fsl_2(\bR)$ or to a Cartan subalgebra of $\fsl_2(\bR)$, otherwise $\gh$ would preserve an isotropic line and therefore be a subalgebra of $\fp_2$. We may also assume that $A_0$ is not contained in a Cartan subalgebra, since $\fh$ has infinite type. 

It follows that
$
A=A_0=\fsl_2(\bR)
$ and $\fh$ is one of the Lie algebras
$$
\fsl_2(\bR)\oplus\fsl_2(\bR)\;,\qquad \fsl_2(\bR)\oplus\bR\mathrm{diag}(1,-1)\;,\qquad\fsl_2(\bR)\oplus\fso_2(\bR)\;.
$$ 
Again $gr(\gg)=V+\fh$, since it is finite-dimensional. We established the first claim of Theorem \ref{thm:dichotomy}, the last claim is an immediate consequence.
\label{sec:dichotomy}
\section{Finite-dimensional    subalgebras  of   $\gp_1^{(\infty)}$  and   $\gp_2^{(\infty)}$, and their associated homogeneous symplectic $4$-manifolds}
\label{sec:K}
In this section, we will consider finite-dimensional transitive Lie algebras of symplectic vector fields (on the symplectic $4$-dimensional space $V$) which are nonlinear and with the isotropy subalgebra of infinite type. Our arguments rely on a closer
look at the full prolongations of the maximal parabolic subalgebras of $\gsp(V)$.

{\it Notation:} For any Lie algebra $\gg$ and associative algebra $A$, we will denote by $A\cdot\gg$ the tensor product $A\otimes\gg$ of $A$ and $\gg$ with its natural structure of Lie algebra.
\subsection{Geometric preliminaries on transitive Lie algebras of vector fields}
\label{sec:preldis}
In this section, we introduce some basic geometric notions on transitive Lie algebras $\gg$ of vector fields on $\bR^n$. (In this paper, we will be interested only in $n\leq 4$.) In our conventions, the origin of $\bR^n$ is always a regular point for $\gg$.

Let $\gg$ be a transitive Lie algebra of (analytic) vector fields on $\bR^n$.
The following notions are a direct consequence of classical facts on imprimitive Lie algebras (see, e.g., \cite[\S 1.5]{SS} and \cite[\S 2]{GKO}).
\begin{Def}
The Lie algebra $\gg$ is {\rm primitive on
an open subset} $\mathcal U\subset \bR^n$ if there
is no foliation of $\mathcal U$ whose leaves are permuted by the (local) one-parameter subgroups corresponding to elements in $\gg$.
Equivalently,  $\gg$ is primitive on $\mathcal U$ if there is no non-trivial involutive distribution $\mathcal D$ on $\mathcal U$ left invariant by $\gg$. Otherwise, $\gg$ is called {\rm imprimitive on} $\mathcal U$.
\end{Def}
\begin{Def}
The Lie algebra $\gg$ is called {\rm primitive} if it is primitive in every open subset $\mathcal U\subset \bR^n$. Otherwise, $\gg$ is {\rm imprimitive}.
\end{Def}
 Let $\gg$ be imprimitive and $x^1,\ldots,x^k$, $y^1,\ldots,y^{n-k}$ be local coordinates such that the
system of imprimitivity corresponding to $\mathcal D$ is given by the leaves $y=const.$ Then any vector field $X\in\gg$ is (locally) of the form
$$
X=f^1\frac{\partial}{\partial x^1}+\cdots+f^k\frac{\partial}{\partial x^k}+g^1\frac{\partial}{\partial y^1}+\cdots+g^{n-k}\frac{\partial}{\partial y^{n-k}}\;,
$$
where the $g$'s are functions of the $y$'s alone. The set  of $X$ with $g^1=\cdots=g^{n-k}=0$ forms an ideal $\mathfrak i$ in $\gg$, which consists of the vector fields leaving all the leaves invariant. In general, it is not a transitive Lie algebra of vector fields when acting  on a leaf.
\begin{Def}
The ideal $\mathfrak i$ is called the {\rm canonical ideal} (of the imprimitive Lie algebra $\gg$, associated to the distribution $\mathcal D$).
\end{Def}
Since the passage to the quotient space of a manifold under the equivalence relation $\thicksim$ induced by a foliation commonly
entails technical difficulties, the following concepts are only locally defined around the origin.
\begin{Def}
\label{def:tranprim}
Let $\gg$ be a transitive Lie algebra of vector fields on $\bR^n$ which is imprimitive on a neighborhood $\mathcal U$ of the origin. Then we say that $\gg$ is:
\begin{enumerate}
\item {\rm transitive on leaves} if the canonical ideal $\mathfrak i$ acting on any leaf $y=const.$ is a transitive Lie algebra;
	\item {\rm primitive on leaves} if it is transitive on leaves and the canonical ideal $\mathfrak i$ acting on any leaf $y=const.$ is a primitive Lie algebra;
	\item {\rm transversally primitive} if the (transitive) Lie algebra $\mathfrak g/\mathfrak i$ acting on the local quotient space $\bR^{n-k}\cong (\mathcal U/\!\thicksim)$ is a primitive Lie algebra.
\end{enumerate}
\end{Def}
\subsection{Finite-dimensional subalgebras of $\gp_{2}^{(\infty)}$}
\label{sec:4.2}
We recall that the Lie algebra
$\gp_{2}^{(\infty)}$ of formal symplectic vector fields is graded and with isotropy subalgebra $\gp_2$.
The parabolic $\gp_2$ is the stabilizer of the line $\bR p_1$, hence it leaves invariant also the $3$-dimensional orthogonal space $p_1^\perp=\bR p_1+ W$, $W=\mathrm{span}(p_2,q_2)$.

It follows that $\gp_2^{(\infty)}$ is of the ``very imprimitive'' sort:
there exists a flag of involutive distributions
\begin{equation}
\label{eq:flagdis}
0\subset \mathcal D_1 \subset \mathcal D_3\subset T\bR^4\;,
\end{equation}
where $\rk\mathcal D_i=i$ for $i=1,3$, which is left invariant by $\fp_2^{(\infty)}$. The systems of imprimitivity corresponding to these distributions are defined on the entire $\bR^4$ and we let $\mathfrak i_1$, $\mathfrak i_3$ be the associated canonical ideals. Clearly $\mathfrak i_1\subset \mathfrak i_3$.

Now, the Lie algebra $\gp_2^{(\infty)}$ is obviously not primitive on the leaves $y=const.$ of $\mathcal D_3$ (here and throughout this section, $y$ is a single real coordinate) and
it is transversally primitive across these leaves in a trivial way.
The action on each leaf of the canonical ideal $\mathfrak i_3$ is transversally primitive (across the leaves of $\mathcal D_1$), at least formally. To check this, we need the following result, a direct consequence of the natural identification of the
algebra of formal symplectic vector fields with the symmetric algebra of $V$ (modulo constants)
and equation \eqref{eq:fullprolpar2} with $p_1=y$ and $q_1=\frac{\partial}{\p y}$.
\begin{Prop}
\label{prop:fullprol2}
The full prolongation of $\mathfrak{p}_2$ is isomorphic, as an abstract Lie algebra, to the semi-direct sum
\begin{equation}
\label{eq:identification2}
\begin{aligned}
\gp_2^{(\infty)}&=\operatorname{Der}(\bR[[y]])+\mathfrak i_3\qquad\text{and}\\
\mathfrak i_3&=
\Big(
\bR[[y]]\otimes \gsp(W)^{(\infty)}+
\bR_+[[y]]\Big)
\end{aligned}
\end{equation}
of the Lie algebra of formal vector fields on the line acting naturally on the canonical ideal $\mathfrak i_3$.
The canonical ideal $\mathfrak i_1$ is the center of $\mathfrak i_3$ and it consists of the formal power series in one variable
modulo constants:
$$
\mathfrak i_1=\bR_+[[y]]\cong\frac{\bR[[y]]}{\bR 1}\;.
$$
Finally
\begin{equation}
\label{eq:contactsympl}
[y^iX,y^jY]=
\begin{cases}
[X,Y]&\text{if}\quad i=j=0,\\
y^{i+j}[X,Y]+y^{i+j}\Omega(X,Y)\quad&\text{if}\quad i+j>0,
\end{cases}
\end{equation}
for all $X,Y\in\gsp(W)^{(\infty)}$, where $\Omega$ is the standard symplectic form at the origin.
\end{Prop}
\bc
The Lie algebra $\mathfrak i_3$ acts transversally primitively on each leaf $y=const.$
\ec
\pf
The (formal) one-parameter family of transverse actions is given by the natural structure of a Lie algebra of $\mathfrak i_3/\mathfrak i_1\cong \bR[[y]]\cdot \gsp(W)^{(\infty)}$ as $y$-dependent symplectic vector fields on $W$.
\qed

The proof of the following lemma is straightforward.

\begin{Lem}
\label{lem:stab2}
The stability subalgebra $\gp_2^{\geq 0}$ (resp. the isotropy subalgebra $\gp_2$) of $\gp_2^{(\infty)}$ correspond, under the identification \eqref{eq:identification2}, to the subalgebras
\begin{equation}
\begin{aligned}
\gp_2^{\geq 0}&= \bR_+[[y]]\frac{\partial}{\p y}+\Big(\bR_+[[y]]\otimes \gsp(W)^{(\infty)}+\gsp(W)^{\geq 0}+\span(y^2,y^3,y^4,\ldots)\Big)\;,\\
&\\
\gp_2&=\bR y \frac{\partial}{\p y}+\Big(yW+\gsp(W)+\bR y^2\Big)\;,
\end{aligned}
\end{equation}
where $\gsp(W)^{\geq 0}$ is the stability subalgebra of $\gsp(W)^{(\infty)}$.
\end{Lem}
We now turn to finite-dimensional transitive subalgebras $\gg$ of $\gp_{2}^{(\infty)}$. As for the entire Lie algebra,
$\gg$ is not primitive on the leaves of the distribution $\mathcal D_3$ and, as a consequence of the next Proposition \ref{prop:vfline}, it is transversally primitive across these leaves. In this context, there is only one meaningful notion of primitivity: we are led to consider the subalgebras $\gg$ whose canonical ideal associated to $\mathcal D_3$ acts transversally primitively (across the leaves of $\mathcal D_1$) on each leaf $y=const.$
\vskip0.2cm\par
We assume, for simplicity, that our finite-dimensional and transitive subalgebra $\gg\subset\gp_2^{(\infty)}$ is consistent with the decomposition \eqref{eq:identification2}, i.e., it is of the form
\begin{equation}
\label{eq:consdec}
\gg=\ga+\widetilde\gg+\xi
\end{equation}
with
\begin{equation}
\label{eq:decompp2}
\ga=\operatorname{Der}(\bR[[y]])\cap\gg\;,\qquad \widetilde\gg=
\Big(\bR[[y]]\otimes \gsp(W)^{(\infty)}\Big)\cap\gg\;,\qquad\xi=\bR_+[[y]]\cap\gg\;.
\end{equation}
We remark that $\widetilde\gg$ is {\it not} a subalgebra of $\gg$ in general, as $[\widetilde\gg,\widetilde\gg]\subset \widetilde\gg+\xi$ due to \eqref{eq:contactsympl}. Nonetheless, we shall often identify $\widetilde\gg$ with the quotient of $\widetilde\gg+\xi$ by its central ideal $\xi$, the Lie brackets being \eqref{eq:contactsympl} with the symplectic form at the origin set to zero, i.e., the Lie brackets of $\bR[[y]]\cdot \gsp(W)^{(\infty)}$.

Transitivity of $\gg$ reads as the transitivity of $\ga$ and $\widetilde\gg$ (in the sense that
the natural projection of $\widetilde\gg$ to $W$ is surjective) and the fact that $\xi$ has a formal power series of the form $y+h.o.t.$
\begin{Def}
\label{def:splitting}
We say that a subalgebra $\gg$ of $\gp_2^{(\infty)}$ is {\rm splitting} if:
\begin{enumerate}
\item it is consistent with the decomposition \eqref{eq:identification2} of $\gp_2^{(\infty)}$, i.e., it is of the form \eqref{eq:consdec};
\item if a vector field $\displaystyle X=\sum_{i=0}^{+\infty}y^iX_i$, $X_i\in\gsp(W)^{(\infty)}$ for all $i\geq 0$, is in $\widetilde\gg$, then its component $X_0$ is in $\widetilde\gg$ too.
\end{enumerate}
\end{Def}
The classification of general (nonsplitting) subalgebras requires studying different deformation and extension problems; in this paper, we will consider only the splitting case. We depart with the following.
\bp
\label{prop:vfline}
(see, e.g., \cite{Dra})
The    following list  exhausts, up to isomorphism,  the  finite-dimensional transitive subalgebras $\ga$ of the Lie algebra of vector  fields  on  the  real line:
\begin{enumerate}
	\item $\ga_1 =  \span(\frac{\partial}{\p y})$,
	\item $\ga_2 = \mathfrak{aff}(\bR) = \span(\frac{\partial}{\p y}, y\frac{\partial}{\p y})$,
	\item $ \ga_3 = \mathfrak{proj}(\bR)= \span(\frac{\partial}{\p y}, y\frac{\partial}{\p y}, y^2\frac{\partial}{\p y})$.
\end{enumerate}
\ep
We recall that an isomorphism of Lie algebras of vector fields on the line is by definition
the map induced by an analytic coordinate change $\phi:\bR[[y]]\to\bR[[y]]$ at the origin. In particular, any such isomorphism acts on $\gp_2^{(\infty)}$ too, in a way compatible with the decomposition into sum  of three Lie algebras. Hence, without any loss of generality, we may take the component $\ga$ of the Lie algebra \eqref{eq:consdec} as in Proposition \ref{prop:vfline} from now on, in particular $\frac{\partial}{\p y}\in\ga$.
\bl
\label{lem:xi}
\hfill\par\noindent
\begin{enumerate}
\item  Any non-zero finite-dimensional subspace  $\xi$ invariant w.r.t. $\ga_2$  has  the  form  $$\xi=P^k_+ = \span(y, y^2, \ldots, y^k)\;.$$ The same holds for $\ga_1$ if, in addition, $\xi$ consists of polynomials;
\item There is no finite-dimensional subspace $\xi$ which is  invariant  w.r.t.  $\ga_3$.
\end{enumerate}
\el
\pf
\begin{enumerate}
	\item Using $y\frac{\partial}{\p y}\in\ga_2$ and Vandermonde determinants, one sees that
any subspace $\xi$ which is $\ga_2$-invariant and does not consist just of polynomials is necessarily infinite-dimensional. The rest follows from $\frac{\partial}{\p y}\in\ga$.
\item Clear. \qed
\end{enumerate}
Lemma \ref{lem:xi} implies that $\ga=\ga_3$ is not permissible as a subalgebra of a splitting finite-dimensional transitive $\gg=\ga+\widetilde\gg+\xi$.
In view of this lemma, we will also enforce one growth condition on $\gg$, namely, we will assume from now on that $\gg$ {\it is polynomial in} $y$.

Let $P^k\subset \bR[[y]]$ be the space of all polynomials of degree $\leq k$. Since $\widetilde\gg$ is finite-dimensional and polynomial in $y$, there is a nonnegative integer $N$ such that
$$\widetilde\gg\subset P^{N}\otimes\gsp(W)^{(\infty)}\;.$$
 We assume that $N$ is the smallest integer with this property and decompose any $X\in\widetilde\gg$ into
\begin{equation}
X=\sum_{i=0}^{N}y^iX_i\;,
\label{eq:finiteX}
\end{equation}
where $X_i\in\gsp(W)^{(\infty)}$ for all $i=0,\ldots, N$. We also set $\overline\gg=\widetilde\gg\cap \gsp(W)^{(\infty)}$ and
$$
\widetilde\gg^i=\Big\{Y\in \gsp(W)^{(\infty)}\mid\text{there exists}\;\;X\in\widetilde\gg\;\;\text{with}\;\;X_i=Y\Big\}
$$
for all $i=0,\ldots, N$. Clearly $\overline\gg\subset \widetilde\gg^0$ and, by the second condition of Definition \ref{def:splitting}, $\overline\gg=\widetilde\gg^0$.
\bp
\label{prop:idealgg}
\hfill\par\noindent
\begin{enumerate}
\item  The space $\widetilde\gg^i$ is contained in $\widetilde\gg^{i-1}$ and it is an ideal of $\overline\gg$, for all $i>0$;
	\item If $\overline{\gg}$ has no non-zero nilpotent ideals then $\widetilde\gg=\overline\gg$.
\end{enumerate}
\ep
\pf
The first claim follows by repeatedly applying $\frac{\partial}{\p y}\in\ga$ to elements \eqref{eq:finiteX} of $\widetilde \gg$.
It remains to show that $\widetilde\gg\subset\overline\gg$ if $\overline{\gg}$ has no non-zero nilpotent ideals. This follows immediately, as the ideal $\widetilde\gg^N$ of $\overline \gg$ is non-zero and abelian if $N>0$.
\qed
The classification of the finite-dimensional Lie algebras of vector fields in the real plane can be found in \cite[Table 1]{GKO} (the notation therein is $p=\partial_x$ and $q=\partial_y$; caution: the origin is not a regular point for the Lie algebras $(2)$ and $(17)$-$(19)$,  a translation of coordinates has to be performed to match our present set-up).

There are transitive Lie algebras of symplectic vector fields, both primitive and not primitive.
The primitive ones are the unimodular affine Lie algebra $\fsl_2(\bR)+\bR^2$ (case (5) in \cite[Table 1]{GKO}), its subalgebra $\bR+\bR^2$ of Euclidean motions (case (1) for $\alpha=0$) and the two filtered deformations $\fsl_2(\bR)$ and
$\gso_3(\bR)$ of the latter given by the Lie algebra of infinitesimal automorphisms $(2)$ of the hyperbolic plane and $(3)$ of the Euclidean $2$-sphere.

The imprimitive ones are the Lie algebras $(17)$-$(18)$, the Lie algebra $(24)$ and its subalgebra $(12)$ for $\alpha=-1$ and, finally, the Lie algebra $(22)$ for an appropriate choice of the defining functions $\eta_i(x)$ --- e.g. $\eta_{i}(x)=x^{i-1}$ for all $i=1,\ldots,r$.

In    all cases, there  is   a natural decomposition
$$\overline{\gg}  = \gs  + \gn$$
into the semidirect  sum of a subalgebra $\gs$ and a (possibly trivial) abelian ideal $\gn$. The latter coincides with the maximal nilpotent ideal of $\overline{\gg}$, unless $\overline{\gg}$ is $(24)$ or $(22)$ with all polynomial defining functions as above.
\vskip0.2cm\par
We now state the main result of this section. Recall the definition of a transversally primitive Lie algebra across the flag of distributions \eqref{eq:flagdis} on $\bR^4$, cf. Definition \ref{def:tranprim} and the discussion after Lemma \ref{lem:stab2}.
\bt
\label{thm:K1}
Let $\overline\gg=\gs+\gn$ be one of the four finite-dimensional primitive Lie algebras of symplectic vector fields on the real plane, as included in \cite[Table 1]{GKO}. Then:
\vskip0.3cm\par\noindent
1. If $\overline\gg=\gs$ is the Lie algebra $\fsl_2(\bR)$ of infinitesimal automorphisms of the hyperbolic plane or the Lie algebra
$\gso_3(\bR)$ of infinitesimal automorphisms of the Euclidean $2$-sphere, then
\begin{equation}
\gg=\mathfrak{aff}(\bR)+\overline\gg+P^k_+
\label{eq:1p2}
\end{equation}
is a finite-dimensional transitive and transversally primitive Lie algebra of symplectic vector fields on the $4$-dimensional space. Here
\begin{equation}
\label{eq:notationaff}
\mathfrak{aff}(\bR) = \span(\frac{\partial}{\p y}, y\frac{\partial}{\p y})\qquad\text{and}\qquad P^k_+ = \span(y, y^2, \ldots, y^k)\;,
\end{equation}
for some integer $k\geq 1$. The associated stability subalgebra $\gk$ and isotropy algebra $\gh$ are given by
\begin{equation}
\label{eq:stabisoI}
\begin{aligned}
\gk&=\begin{cases}
\bR y\frac{\partial}{\p y}+\Big(\overline\gk+\span(y^2,\ldots,y^k)\Big)\qquad&\text{if}\;\;k>1\\
&\\
\bR y\frac{\partial}{\p y}+\overline\gk\qquad&\text{if}\;\;k=1
\end{cases}
&\\
&\\
\gh&=\begin{cases}
\bR y \frac{\partial}{\p y}+\Big(\overline\gh+\bR y^2\Big)\qquad&\text{if}\;\;k>1\\
&\\
\bR y \frac{\partial}{\p y}+\overline\gh\qquad&\text{if}\;\;k=1
\end{cases}\quad\;,
\end{aligned}
\end{equation}
where $\overline\gk$ is the stability subalgebra of $\overline\gg$ and $\overline\gh\cong \gso_2(\bR)$ the linear isotropy algebra of $\overline\gg$.
\vskip0.3cm\par\noindent
2. If $\overline\gg=\gs+\gn$ is the unimodular affine Lie algebra $\fsl_2(\bR)+\bR^2$ or its subalgebra $\bR+\bR^2$ of Euclidean motions, then the nilradical $\gn=\bR^2$ of $\overline\gg$ is abelian and
\begin{equation}
\label{eq:2p2}
\gg=\mathfrak{aff}(\bR)+(\gs+P^N\otimes\gn)+P^k_+
\end{equation}
is a finite-dimensional transitive and transversally primitive Lie algebra of symplectic vector fields on the $4$-dimensional space, where the notation is as in \eqref{eq:notationaff} and
$$P^N=\span(1,y, y^2, \ldots, y^N)$$
for some nonnegative integer $N$ with $2N\leq k$. The stability subalgebra $\gk$ is given by
\begin{equation}
\label{eq:uffinalI}
\begin{aligned}
\gk&=
\begin{cases}
\bR y\frac{\partial}{\p y}+\Big(\gs+P^N_+\otimes \gn+\span(y^2,\ldots,y^k)\Big)\qquad&\text{if}\;\;k>1\;\;\text{and}\;\;N>0\\
&\\
\bR y\frac{\partial}{\p y}+\Big(\gs+\span(y^2,\ldots,y^k)\Big)\qquad&\text{if}\;\;k>1\;\;\text{and}\;\;N=0\\
&\\
\bR y\frac{\partial}{\p y}+\gs\qquad&\text{if}\;\;k=1
\end{cases}
\end{aligned}
\end{equation}
and the isotropy algebra $\gh$ by
\begin{equation}
\label{eq:uffinalII}
\begin{aligned}
\gh&=
\begin{cases}
\bR y \frac{\partial}{\p y}+\Big(\gs+y\gn+\bR y^2\Big)\qquad&\text{if}\;\;k>1\;\;\text{and}\;\;N>0\\
&\\
\bR y \frac{\partial}{\p y}+\Big(\gs+\bR y^2\Big)\qquad&\text{if}\;\;k>1\;\;\text{and}\;\;N=0\\
&\\
\bR y \frac{\partial}{\p y}+\gs\qquad&\text{if}\;\;k=1
\end{cases}
\end{aligned}
\end{equation}
The Lie algebras \eqref{eq:1p2} and \eqref{eq:2p2} just described are subalgebras of $\gp_2^{(\infty)}$, the explicit expression of their Lie brackets is as in Proposition \ref{prop:fullprol2}.
\vskip0.3cm\par
Conversely:
\vskip0.3cm\par\noindent
3. Any finite-dimensional transitive and transversally primitive subalgebra $\gg$
of $\gp_2^{(\infty)}$ that is splitting and polynomial in $y$ is a subalgebra
\begin{equation}
\label{eq:useat4}
\gg=\ga+\widetilde\gg+\xi
\end{equation}
of \eqref{eq:1p2} or \eqref{eq:2p2}. More precisely, $\ga=\span(\frac{\partial}{\p y})$ or $\ga=\mathfrak{aff}(\bR)$, $\widetilde\gg\subset\gs+P^N\otimes\gn$ for one of the finite-dimensional primitive Lie algebras $\overline\gg=\gs+\gn$ of symplectic vector fields on the real plane
and $\xi=P^k_+$ for some $k\geq 1$.
\vskip0.3cm\par\noindent
4. For all Lie algebras \eqref{eq:useat4} of symplectic vector fields, there exists a choice of a connected Lie group $G$ with Lie algebra $Lie(G)=\gg$
so that the analytic subgroup $K$ of $G$ with Lie algebra $Lie(K)=\gk=\gg\cap\gp_2^{\geq 0}$ is closed and
$$M = (G/K,\o)$$
is a homogeneous symplectic $4$-manifold. If $k>2$, then $M=G/K$ does not admit any invariant linear connection, in particular it is not reductive and it is not a homogeneous Fedosov manifold.
\et
\pf
The first two claims follow from direct computations using the Lie brackets described in Proposition \ref{prop:fullprol2};
the fact that the given Lie algebras \eqref{eq:1p2} and \eqref{eq:2p2} are transitive and transversally primitive is also immediate. We now show the converse implication.

Let $\gg=\ga+\widetilde\gg+\xi$ be a splitting
subalgebra of $\gp_2^{(\infty)}$ with the desired properties. By transitivity of $\gg$ and the fact that $\gg$ is polynomial in $y$ we have $\ga\subset
\mathfrak{aff}(\bR)$ and $\xi=P^k_+$ for some $k\geq 1$, see Proposition \ref{prop:vfline} and
Lemma \ref{lem:xi}. We now focus on $\widetilde\gg$.

First $\overline\gg=\widetilde\gg\cap \gsp(W)^{(\infty)}$ is a transitive Lie algebra of symplectic vector fields on the real plane,
since the natural projection of $\widetilde\gg$ to $W$ is surjective and the component $X_0$ of any element $\displaystyle X=\sum_{i=0}^{+\infty}y^iX_i$ of $\widetilde\gg$ is still in $\widetilde\gg$. By assumption, $\gg$ is transversally primitive on each leaf $y=const.$ and so, in particular, on the leaf $y=0$. This says that $\overline\gg$ is primitive too, hence one of the four finite-dimensional primitive Lie algebras of symplectic vector fields on the real plane.

If $\overline\gg$ is simple, then $\widetilde\gg=\overline\gg$ by Proposition
\ref{prop:idealgg} and $\gg$ is a subalgebra of \eqref{eq:1p2}. Otherwise $\overline\gg=\fsl_2(\bR)+\bR^2$ or $\overline\gg=\bR+\bR^2$ and we now deal with the two cases separately.

Recall first that
$$\widetilde\gg\subset P^{N}\otimes\gsp(W)^{(\infty)}\;,$$
for some nonnegative integer $N$. The case $N=0$ is trivial, as it amounts to $\widetilde\gg=\overline\gg$. If $N>0$, then by Proposition \ref{prop:idealgg} and the fact that $\widetilde\gg^N$ is a non-zero abelian ideal of $\overline\gg$, we have
$$
\widetilde\gg^i=\begin{cases}
\overline\gg\qquad&\text{for all}\qquad 0\leq i\leq M\\
\gn\qquad&\text{for all}\qquad M+1\leq i \leq N\\
0\qquad&\text{for all}\qquad N<i
\end{cases}
$$
for some nonegative integer $M<N$. We want to show $M=0$, so that $\widetilde\gg\subset \gs+P^N\otimes\gn$.
\begin{itemize}
	\item[(i)] Let $\overline\gg=\gs+\gn=\fsl_2(\bR)+\bR^2$ and assume $M>0$. We take elements
	\begin{equation}
		X=\sum_{i=0}^{N}y^iX_i\qquad\text{and}\qquad Y=\sum_{i=0}^{N}y^iY_i
	\label{eq:elementsXY}
	\end{equation}
of $\widetilde\gg$ such that $X_M, Y_M\in \gs$. Then the component of the bracket
$$[X,Y]_{2M}\equiv [X_M,Y_M] \mod\gn\;,$$
since $\gn$ is an ideal. On the other hand $2M>M$, whence $[X,Y]_{2M}\in\gn$, which is not possible for all $X_M,Y_M\in \gs$.
\item[(ii)] Let $\overline\gg=\gs+\gn=\bR+\bR^2$, and let $J$ be the complex structure generating $\gs$. Assume $M>0$ and take
$X$ and $Y$ in $\widetilde\gg$ such that $X_M=J$ and $Y_N$ is a nonzero element of $\gn$. By possibly replacing $Y$ with $[J,Y]$, we may assume that {\it all} components of $Y$ are in $\gn$. Hence
\begin{equation}
\begin{aligned}
0&=[X,Y]_{N+M}=[X_M,Y_N]\\
&=[J,Y_N]\neq 0\;,
\end{aligned}
\end{equation}
since $N+M>N$ and $\gn$ is an abelian ideal, an absurd.
\end{itemize}
This proves claim 3. of the theorem.
\vskip0.2cm\par
We now turn to the existence of homogeneous symplectic $4$-manifolds. This amounts to show that the putative stability subgroup is closed. First we note that it is sufficient to establish claim 4. for the ``maximal'' Lie algebras \eqref{eq:1p2} and \eqref{eq:2p2} (for the nonmaximal Lie algebras, one simply considers the orbits under smaller group actions).

Let $G$ be the simply connected Lie group with Lie algebra $\gg$ and $K$ the analytic subgroup of $G$ with $Lie(K)=\gk=\gg\cap\gp_2^{\geq 0}$.
Let $K^{Mal}$ be the topological closure of $K$ in $G$ and $Lie(K^{Mal})=\gk^{Mal}$ the corresponding subalgebra of $\gg$. It is called the {\it Malcev closure} of $\gk$ and it satisfies the following fundamental properties (see \cite{GOVI}):
\begin{itemize}
\item $\gk^{Mal}\supset\gk$,
\item $[\gk^{Mal},\gk^{Mal}]=[\gk,\gk]$.
\end{itemize}
We will now show that $\gk^{Mal}=\gk$ in our case. If $\overline\gg=\gs$ is simple, then
$$
[\gk,\gk]=\begin{cases}
\span(y^2,\ldots,y^k)\qquad&\text{if}\;\;k>1\\
&\\
0\qquad&\text{if}\;\;k=1
\end{cases}
$$
as the stability subalgebra $\overline\gk$ of $\overline\gg$ is $1$-dimensional compact. The subalgebra $\gk$ has codimension four in $\gg$ and
$$
\gg/\gk\cong \span\Big(\frac{\partial}{\p y}, y\Big)+\overline\gg/\overline\gk
$$
as vector spaces. Let $X\in\gk^{Mal}$ be of the form
$$X=\lambda\frac{\partial}{\p y}+\mu y+Z\;,$$
where $Z\in \overline\gg/\overline\gk$. Then the bracket
$$
[X,y\frac{\partial}{\p y}]=\lambda\frac{\partial}{\p y}-\mu y
$$
belongs to $[\gk,\gk]$, since $y\frac{\partial}{\p y}\in\gk$, but this is possible only if $\lambda=\mu=0$. Now, the
hyperbolic plane and the Euclidean $2$-sphere are irreducible Riemannian symmetric spaces and the stability subalgebra $\overline\gk\cong \fso_2(\bR)$ is a maximal subalgebra of $\overline\gg$. Hence either $\gk^{Mal}=\gk$ or
$$
\gk^{Mal}=\begin{cases}
\bR y\frac{\partial}{\p y}+\Big(\overline\gg+\span(y^2,\ldots,y^k)\Big)\qquad&\text{if}\;\;k>1\\
&\\
\bR y\frac{\partial}{\p y}+\overline\gg\qquad&\text{if}\;\;k=1
\end{cases}
$$
but in the second case $[\gk,\gk]\supset\overline\gg$, which is not possible. Hence $\gk^{Mal}=\gk$ and $K=K^{Mal}$ since
both Lie groups are connected. The case $\overline\gg=\gs+\gn$ not simple is similar and we omit the proof.

Finally we use \eqref{eq:stabisoI} and \eqref{eq:uffinalI}-\eqref{eq:uffinalII} to see that the isotropy representation is not exact if $k>2$ so that the associated homogeneous manifolds do not admit invariant linear connections.
\qed
An immediate consequence of Theorem \ref{thm:K1} is the following.
\bc
There exist finite-dimensional transitive subalgebras $\gg$ of $\gp_2^{(\infty)}$ of any positive order and with associated isotropy algebra of infinite type.
\ec
\pf
The Lie algebras \eqref{eq:1p2} and \eqref{eq:2p2} with any $k>1$ do the job. It is not difficult to see that the isotropy subalgebra includes a rank one element and it is therefore of infinite type.
\qed
\subsection{Finite-dimensional subalgebras of $\gp_{1}^{(\infty)}$}
\label{sec:4.3}
The isotropy algebra $\gp_1$ of $\gp_{1}^{(\infty)}$ is the stabilizer
of the Lagrangian plane $P= \mathrm{span}(p_1,p_2)$. Hence,
there exists an involutive Lagrangian distribution
\begin{equation}
\label{eq:flagdisII}
0\subset \mathcal D\subset T\bR^4
\end{equation}
of $\rk\mathcal D=2$ which is left invariant by $\fp_1^{(\infty)}$. The system of imprimitivity of this distribution is defined on the entire $\bR^4$ and we denote by $\mathfrak i$ the associated canonical ideal (see \S \ref{sec:preldis}).
\vskip0.2cm\par
The following is an alternative characterization of $\gp_1^{(\infty)}$.
\begin{Prop}
\label{prop:fullprol1}
The full prolongation of $\mathfrak{p}_1$ is isomorphic, as an abstract Lie algebra, to
the semi-direct sum
\begin{equation}
\label{eq:abstractp1}
\begin{aligned}
\gp_{1}^{(\infty)}&=\operatorname{Der}(\bR[[x,y]])+\mathfrak i
\qquad\text{and}\\
\mathfrak i&=
\bR_+[[x,y]]\cong \frac{\bR[[x,y]]}{\bR 1}
\end{aligned}
\end{equation}
of the Lie algebra of formal vector fields on the plane acting naturally on the canonical ideal. The latter is abelian and identifiable with the formal power series in two variables, modulo constants.

The stability subalgebra $\gp_1^{\geq 0}$ (resp. the isotropy subalgebra $\gp_1$) of $\gp_1^{(\infty)}$ correspond, under this identification, to the subalgebras
\begin{equation}
\begin{aligned}
\gp_1^{\geq 0}&=\operatorname{Der}^{\geq 0}(\bR[[x,y]])+\span(x^2,xy,y^2,\ldots)\;,\\
&\\
\gp_1&=\ggl_2(\bR)+\span(x^2,xy,y^2)\;,
\end{aligned}
\end{equation}
where $\operatorname{Der}^{\geq 0}(\bR[[x,y]])$ is the stability subalgebra of $\operatorname{Der}(\bR[[x,y]])$.
\end{Prop}
In particular the Lie algebra $\gp_1^{(\infty)}$ acts in a transversally primitive fashion across the leaves of $\mathcal D$ (cf. Definition \ref{def:tranprim}; in our case, this means that the quotient $\gg/\mathfrak i$ acts primitively on the plane).

Let $\gg$ be a finite-dimensional subalgebra of $\gp_1^{(\infty)}$ consistent with the above decomposition, that is
$$
\gg=\widetilde\gg+\xi
$$
with
\begin{equation}
\label{eq:decompp1}
\widetilde\gg=
\operatorname{Der}(\bR[[x,y]])\cap\gg\;,\qquad\xi=\bR_+[[x,y]]\cap\gg\;.
\end{equation}
We call any such $\gg$ {\it splitting}. The classification of nonsplitting subalgebras of $\gp_1^{(\infty)}$ requires the study of a
deformation problem; here, we will consider only splitting subalgebras.
\begin{Prop}
Let $\widetilde\gg$ be any transitive Lie algebra of vector fields on the real plane and $\xi$ a subspace of $\bR_+[[x,y]]$ which is $\widetilde\gg$-invariant and such that $x+h.o.t.$ and $y+h.o.t.$ are in $\xi$. Then $$\gg=\widetilde\gg+\xi$$ is a transitive splitting subalgebra of $\gp_1^{(\infty)}$. Conversely, any transitive splitting subalgebra of $\gp_1^{(\infty)}$ is obtained in this way.
The Lie algebra $\gg$ acts in a transversally primitive way across the leaves of $\mathcal D$ if and only if $\widetilde\gg$ acts primitively on the plane.
\end{Prop}
The list of finite-dimensional primitive Lie algebras of vector fields on the real plane is in \cite[Table 1, I]{GKO}. (We emphasize once again that in \cite{GKO} the origin is not a regular point for the Lie algebra of infinitesimal automorphisms of the hyperbolic plane). They all consist of polynomial vector fields. For simplicity,  we also assume that $\xi$ consists of polynomials.
\begin{Lem}
\hfill\par\noindent
\begin{enumerate}
\item
If $\widetilde\gg$ is the Lie algebra of (unimodular) affine infinitesimal transformations of the plane, then any nontrivial finite-dimensional polynomial $\widetilde\gg$-module $\xi$ has the form
$$\xi=P^k_+ = \span(x,y,x^2,xy,y^2,\ldots, x^k,\ldots,y^k)\;;$$
\item
Let $\widetilde\gg$ be the Lie algebra $\gsl_2(\bR)$ of infinitesimal automorphisms of the hyperbolic plane,
$\fso_3(\bR)$ of the Euclidean $2$-sphere, $\fsl_3(\bR)$ of the projective plane and $\fso(1,3)$ of the conformal $2$-sphere. Then
any finite-dimensional polynomial $\widetilde\gg$-module $\xi$ is trivial, i.e., $\xi=0$.
\end{enumerate}
\end{Lem}
\pf
The first claim is immediate and we omit the proof. Note also that it is sufficient to establish the second claim for the Lie algebras $\widetilde\gg=\fsl_2(\bR)$, $\fso_3(\bR)$ and $\fsl_3(\bR)$.

Let $\xi\neq 0$ and take a nonconstant polynomial $f$ of maximum degree in $\xi$, say $\deg(f)=k$, which we decompose into the sum $f=\sum_{i=1}^kf_i$ of (nonconstant) homogeneous polynomials.
If $\widetilde\gg=\fsl_2(\bR)$ or $\fso_3(\bR)$, we may act on such polynomial with a vector field of the form
$$X=(x^2-y^2)\frac{\partial}{\p x}+2xy\frac{\partial}{\p y}+l.o.t.$$
The resulting component  of highest degree is given by $$(x^2-y^2)\frac{\partial}{\p x}(f_k)+2xy\frac{\partial}{\p y}(f_k)$$ and its vanishing implies $f_k=0$, an absurd. The case $\widetilde\gg=\fsl_3(\bR)$ is similar.
\qed
It remains to give $\xi$ for the Lie algebra of infinitesimal conformal
automorphisms of the plane and a $1$-parameter family of deformations of the subalgebra of Euclidean motions. We set
$$E=x\frac{\partial}{\p x}+y\frac{\partial}{\p y}\qquad\text{and}\qquad J=x\frac{\partial}{\p y}-y\frac{\partial}{\p x}\;,$$
so that
\begin{itemize}
\item[(i)] $\widetilde\gg=\mathfrak{conf}(\bR^2)=\mathrm{span}(\frac{\partial}{\p x},\frac{\partial}{\p y},E,J)$,
\item[(ii)] $\widetilde\gg=\mathfrak{euc}_\alpha(\bR^2)=\mathrm{span}(\frac{\partial}{\p x},\frac{\partial}{\p y},J_{\alpha}=\alpha E-J)$,
\end{itemize}
for any $\alpha\geq 0$.
It is clear that the action of $\bC=\mathrm{span}(E,J)$ on the complexified formal power series is diagonal, with $1$-dimensional eigenspaces
$$P^{k,\ell}=\mathrm{span}(z^{\tfrac{k+\ell}{2}}\cdot\overline z^{\tfrac{k-\ell}{2}})\;,\qquad E|_{P^{k,\ell}}=k\Id\quad\text{and}\quad J|_{P^{k,\ell}}=i\ell\Id\;.$$
Here we are assuming $\ell=k,k-2,\ldots,2-k,-k$, otherwise we set $P^{k,\ell}=0$. In this way, the space of complexified formal power series inherits a $\bZ$-bigrading of the form
$
\bC[[x,y]]=\bigoplus_{k,\ell} P^{k,\ell}
$.
Consider the ``upward triangle'' with bigraded nodes the $1$-dimensional eigenspaces and top node given by $P^{k_o,\ell_o}$, pictorially represented by
$$
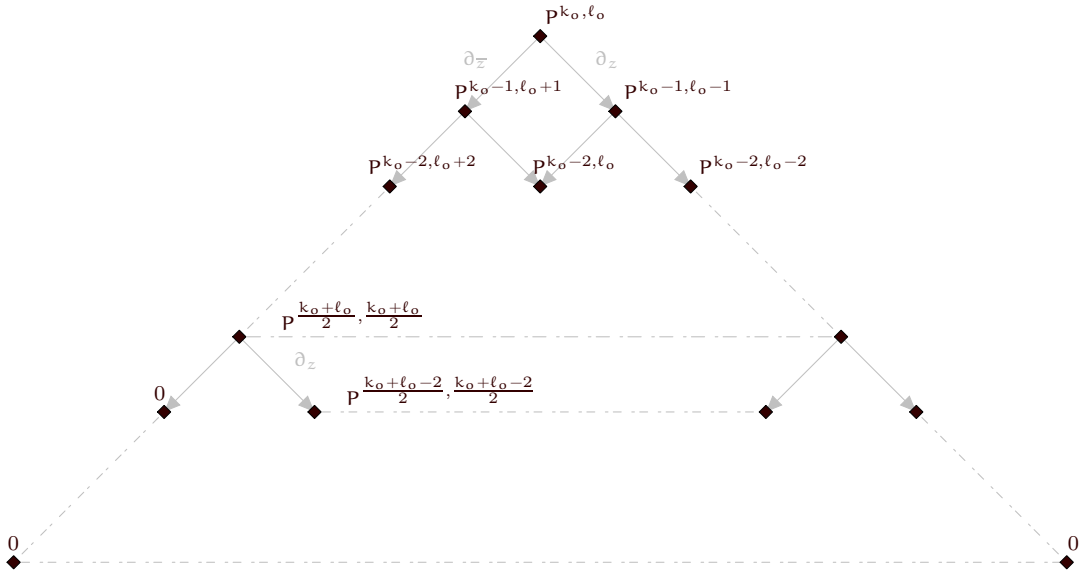
\begin{tikzpicture}[line cap=round,line join=round,>=triangle 45,x=1.0cm,y=1.0cm]
\clip(-7.645547721739148,-2.41795670103093) rectangle (7.632090330434772,5.431798625429551);
\draw [line width=0.4pt,dash pattern=on 1pt off 2pt on 3pt off 4pt,color=cqcqcq] (-2.,3.)-- (-4.,1.);
\draw [line width=0.4pt,dash pattern=on 1pt off 2pt on 3pt off 4pt,color=cqcqcq] (2.,3.)-- (4.,1.);
\draw [line width=0.4pt,dash pattern=on 1pt off 2pt on 3pt off 4pt,color=cqcqcq] (5.,0.)-- (7.,-2.);
\draw [line width=0.4pt,dash pattern=on 1pt off 2pt on 3pt off 4pt,color=cqcqcq] (-7.,-2.)-- (7.,-2.);
\draw [line width=0.4pt,dash pattern=on 1pt off 2pt on 5pt off 4pt,color=cqcqcq] (-4.,1.)-- (4.,1.);
\draw [->,line width=0.4pt,color=cqcqcq] (-1.,4.) -- (0.,3.);
\draw [line width=0.4pt,dash pattern=on 1pt off 2pt on 3pt off 4pt,color=cqcqcq] (-5.,0.)-- (-7.,-2.);
\draw [line width=0.4pt,dash pattern=on 1pt off 2pt on 3pt off 4pt,color=cqcqcq] (-3.,0.)-- (3.,0.);
\draw [->,line width=0.4pt,color=cqcqcq] (-1.,4.) -- (-2.,3.);
\draw [->,line width=0.4pt,color=cqcqcq] (1.,4.) -- (0.,3.);
\draw [->,line width=0.4pt,color=cqcqcq] (1.,4.) -- (2.,3.);
\draw [->,line width=0.4pt,color=cqcqcq] (-4.,1.) -- (-5.,0.);
\draw [->,line width=0.4pt,color=cqcqcq] (-4.,1.) -- (-3.,0.);
\draw [->,line width=0.4pt,color=cqcqcq] (4.,1.) -- (5.,0.);
\draw [->,line width=0.4pt,color=cqcqcq] (0.,5.) -- (-1.,4.);
\draw [->,line width=0.4pt,color=cqcqcq] (0.,5.) -- (1.,4.);
\draw [->,line width=0.4pt,color=cqcqcq] (4.,1.) -- (3.,0.);
\begin{scriptsize}
\draw [fill=ttqqqq] (0.,5.) ++(-2.5pt,0 pt) -- ++(2.5pt,2.5pt)--++(2.5pt,-2.5pt)--++(-2.5pt,-2.5pt)--++(-2.5pt,2.5pt);
\draw[color=ttqqqq] (0.47817025122115814,5.283435552249027) node {$P^{k_o,\ell_o}$};
\draw [fill=ttqqqq] (-1.,4.) ++(-2.5pt,0 pt) -- ++(2.5pt,2.5pt)--++(2.5pt,-2.5pt)--++(-2.5pt,-2.5pt)--++(-2.5pt,2.5pt);
\draw[color=ttqqqq] (-0.420635186389557,4.285356696307316) node {$P^{k_{o}-1,\ell_{o}+1}$};
\draw [fill=ttqqqq] (1.,4.) ++(-2.5pt,0 pt) -- ++(2.5pt,2.5pt)--++(2.5pt,-2.5pt)--++(-2.5pt,-2.5pt)--++(-2.5pt,2.5pt);
\draw[color=ttqqqq] (1.8339567821247503,4.285356696307316) node {$P^{k_{o}-1,\ell_{o}-1}$};
\draw [fill=ttqqqq] (0.,3.) ++(-2.5pt,0 pt) -- ++(2.5pt,2.5pt)--++(2.5pt,-2.5pt)--++(-2.5pt,-2.5pt)--++(-2.5pt,2.5pt);
\draw[color=ttqqqq] (0.4613098618374152,3.2872778403656056) node {$P^{k_{o}-2,\ell_o}$};
\draw [fill=ttqqqq] (-2.,3.) ++(-2.5pt,0 pt) -- ++(2.5pt,2.5pt)--++(2.5pt,-2.5pt)--++(-2.5pt,-2.5pt)--++(-2.5pt,2.5pt);
\draw[color=ttqqqq] (-1.5684312176937765,3.2872778403656056) node {$P^{k_{o}-2,\ell_{o}+2}$};
\draw [fill=ttqqqq] (2.,3.) ++(-2.5pt,0 pt) -- ++(2.5pt,2.5pt)--++(2.5pt,-2.5pt)--++(-2.5pt,-2.5pt)--++(-2.5pt,2.5pt);
\draw[color=ttqqqq] (2.830324481179571,3.2872778403656056) node {$P^{k_{o}-2,\ell_{o}-2}$};
\draw [fill=ttqqqq] (-5.,0.) ++(-2.5pt,0 pt) -- ++(2.5pt,2.5pt)--++(2.5pt,-2.5pt)--++(-2.5pt,-2.5pt)--++(-2.5pt,2.5pt);
\draw[color=ttqqqq] (-5.0588577750019065,0.245834840164852) node {0};
\draw [fill=ttqqqq] (-3.,0.) ++(-2.5pt,0 pt) -- ++(2.5pt,2.5pt)--++(2.5pt,-2.5pt)--++(-2.5pt,-2.5pt)--++(-2.5pt,2.5pt);
\draw[color=ttqqqq] (-1.3152916070775196,0.279553720433153) node {$P^{\tfrac{k_o+\ell_o-2}{2},\tfrac{k_o+\ell_o-2}{2}}$};
\draw [fill=ttqqqq] (3.,0.) ++(-2.5pt,0 pt) -- ++(2.5pt,2.5pt)--++(2.5pt,-2.5pt)--++(-2.5pt,-2.5pt)--++(-2.5pt,2.5pt);
\draw [fill=ttqqqq] (5.,0.) ++(-2.5pt,0 pt) -- ++(2.5pt,2.5pt)--++(2.5pt,-2.5pt)--++(-2.5pt,-2.5pt)--++(-2.5pt,2.5pt);
\draw [fill=ttqqqq] (-4.,1.) ++(-2.5pt,0 pt) -- ++(2.5pt,2.5pt)--++(2.5pt,-2.5pt)--++(-2.5pt,-2.5pt)--++(-2.5pt,2.5pt);
\draw[color=ttqqqq] (-2.4843630406351758,1.2776325763748637) node {$P^{\tfrac{k_o+\ell_o}{2},\tfrac{k_o+\ell_o}{2}}$};
\draw [fill=ttqqqq] (4.,1.) ++(-2.5pt,0 pt) -- ++(2.5pt,2.5pt)--++(2.5pt,-2.5pt)--++(-2.5pt,-2.5pt)--++(-2.5pt,2.5pt);
\draw [fill=ttqqqq] (7.,-2.) ++(-2.5pt,0 pt) -- ++(2.5pt,2.5pt)--++(2.5pt,-2.5pt)--++(-2.5pt,-2.5pt)--++(-2.5pt,2.5pt);
\draw[color=ttqqqq] (7.087409321618137,-1.7503228717185695) node {0};
\draw [fill=ttqqqq] (-7.,-2.) ++(-2.5pt,0 pt) -- ++(2.5pt,2.5pt)--++(2.5pt,-2.5pt)--++(-2.5pt,-2.5pt)--++(-2.5pt,2.5pt);
\draw[color=ttqqqq] (-7.001593173111548,-1.7503228717185695) node {0};
\draw[color=cqcqcq] (-3.1220383680302612,0.6841802836527654) node {$\partial_z$};
\draw[color=cqcqcq] (-0.8767173308430865,4.689983259526929) node {$\partial_{\overline z}$};
\draw[color=cqcqcq] (0.8767173308430865,4.689983259526929) node {$\partial_{z}$};
\end{scriptsize}
\end{tikzpicture}
$$
The partial derivatives
$$
\frac{\partial}{\p z}:P^{k,\ell}\rightarrow P^{k-1,\ell-1}\;,\qquad \frac{\partial}{\p \overline z}:P^{k,\ell}\rightarrow P^{k-1,\ell+1}\;,
$$
sends nodes into nodes, moving downwards right and left, respectively.
We also note that $E$ (resp. $J$) acts with the same eigenvalue on all the nodes on a same vertical (resp. horizontal) level. We denote the direct sum of all the nodes $P^{k,\ell}$ with $k>0$ included in the triangle with top node $P^{k_o,\ell_o}$ by $W^{k_o,\ell_o}$, it is
a (complex) module for $\widetilde\gg=\mathfrak{conf}(\bR^2)$.
\\

The modules $W^{k,\ell}$ just constructed are all inequivalent. This implies the following.
\begin{Prop}
\label{prop:triangle}
\hfill\par\noindent
\begin{enumerate}
\item
Any nontrivial finite-dimensional complex polynomial $\mathfrak{conf}(\bR^2)$-module $\Xi$ is the sum (not necessarily direct)
of a finite number of modules $W^{k,\ell}$. It is invariant by conjugation if and only if it includes both modules $W^{k,\pm\ell}$ or none of them;
\item
Any nontrivial finite-dimensional (real) polynomial $\mathfrak{conf}(\bR^2)$-module $\xi$ is the real form of a unique complex module $\Xi$
invariant by conjugation.
\end{enumerate}
\end{Prop}
Note that $J_\alpha$ acts on $P^{k,\ell}$ with eigenvalue $\alpha k-i\ell$. The modules $P^{k,\ell}$
are still all inequivalent under this action, unless $\alpha=0$ so that $J_\alpha=-J$ and $P^{k,\ell}\cong P^{k',\ell}$.
In any case, Proposition \ref{prop:triangle} carries over with minor modifications and it is not difficult to see that any $\mathfrak{euc}_{\alpha}(\bR^2)$-module is always included in a $\mathfrak{conf}(\bR^2)$-module as we just described.
\vskip0.2cm\par
We have proved most of the following.
 \bt
\label{thm:K2}
Let $\widetilde\gg$ be one of the following finite-dimensional primitive Lie algebras of vector fields on the real plane:
\begin{enumerate}
\item  the Lie algebra $\widetilde\gg=\fsl_2(\bR)+\bR^2$ of unimodular affine infinitesimal transformations;
\item  the Lie algebra $\widetilde\gg=\ggl_2(\bR)+\bR^2$ of affine infinitesimal transformations;
\item  the Lie algebra $\widetilde\gg=\mathfrak{conf}(\bR^2)$ of infinitesimal conformal
automorphisms of the plane;
\item a $1$-parameter family $\widetilde\gg=\mathfrak{euc}_\alpha(\bR^2)$ of deformations of the subalgebra of Euclidean motions.
\end{enumerate}
Then the semi-direct sum
\begin{equation}
\label{eq:quasifine}
\gg=\widetilde\gg+\xi
\end{equation} is a finite-dimensional transitive and transversally primitive subalgebra of $\gp_1^{(\infty)}$, where in the first two cases $\xi=P^k_+$, $k\geq 1$, is the space of nonconstant polynomials in two variables of degree at most $k$  and in the last two cases $\xi$ is a sum of ``triangle'' modules $W^{k,\ell}$ as in Proposition \ref{prop:triangle}. The stability subalgebra $\gk$ and the isotropy algebra $\gh$ are given by
\begin{equation}
\begin{aligned}
\gk&=\widetilde\gk+(\xi\cap\span(x^2,xy,y^2,\ldots))\;,\\
&\\
\gh&=\widetilde\gk+(\xi\cap\span(x^2,xy,y^2))\;,
\end{aligned}
\end{equation}
where $\widetilde\gk$ is the stability subalgebra of $\widetilde\gg$.

Let $\widetilde G=\widetilde K\ltimes\bR^2$ be the connected Lie group of affine transformations of $\bR^2$, where $Lie(\widetilde G)=\widetilde\gg$ and $Lie(\widetilde K)=\widetilde\gk$. Then $K=\widetilde K\ltimes (\xi\cap\span(x^2,xy,y^2,\ldots))$ is a closed subgroup of
$G=\widetilde G\ltimes \xi$ and
$M = (G/K,\o)$ a homogeneous symplectic $4$-manifold. If $k>2$, then $M=G/K$ is not reductive and it is not a homogeneous Fedosov manifold.

Conversely, any finite-dimensional transitive and transversally primitive subalgebra $\gg$
of $\gp_1^{(\infty)}$ that is splitting and polynomial in $x$ and $y$ is a subalgebra of some of the Lie algebras \eqref{eq:quasifine}.
Actually, it coincides with one of them, unless $\widetilde\gg=\mathfrak{euc}_\alpha(\bR^2)$ and $\a=0$, in which case $\xi$ may be properly contained in a sum of ``triangle'' modules.
\et
\pf
It remains to note that the isotropy representation is not exact if $k>2$ and hence $M=G/K$ does not admit any invariant linear connection.
\qed
\bc
There exist finite-dimensional transitive subalgebras $\gg$ of $\gp_1^{(\infty)}$ of any positive order and with associated isotropy algebra of infinite type.
\ec
\section{Homogeneous  torsion free  symplectic  connections  (Fedosov  structures)}
\subsection{Existence of homogeneous Fedosov structures}
\label{sec:lastsecI}
It is a well known fact that on a $2n$-dimensional symplectic manifold $(M,\o)$ there always exists a linear connection $\nabla$ which is torsion free and preserves $\o$. In other words,  the triple $(M,\o,\nabla)$ is a Fedosov manifold, in the terminology of \cite{GRS}.
The usual way to see this is to take trivial local connections in Darboux coordinates and globally glue them using a partition of unity.
The following alternative argument is borrowed from \cite{BCGRS}.

Take $\nabla^o$ {\it any} torsion free linear connection (e.g., the Levi-Civita connection associated to a Riemannian metric on $M$). The covariant derivative $\nabla^o\o$ is a section of $T^*M\otimes \Lambda^2 T^*M$ and, since $\o$ is closed, we have
$$
\mathfrak{S}_{XYZ}\Big\{\nabla_X^o\o(Y,Z)\Big\}=0\;,
$$
where $\mathfrak{S}_{XYZ}$ is the cyclic sum over the vector field $X,Y,Z$ on $M$. We define the section $N$ of $T^*M\otimes T^*M\otimes TM$ by
$
\omega(N(X,Y),Z)=\nabla_X^o\o(Y,Z)
$
and a new connection by
$$
\nabla_XY=\nabla^o_XY+\frac{1}{3}N(X,Y)+\frac{1}{3}N(Y,X)\;.
$$
\vskip-0.6cm\par\noindent
\bp\cite{BCGRS}
\label{prop:fedosovcon}
The connection $\nabla$ is torsion free and preserves $\omega$.
\ep
\pf
It is clear that $\nabla$ is torsion free. We then compute
$$
\begin{aligned}
\nabla_X\o(Y,Z)&=X(\o(Y,Z))-\o(\nabla_X Y,Z)-\o(Y,\nabla_X Z)\\
&=\nabla^o_X\o(Y,Z)-\frac{1}{3}\o(N(X,Y),Z)-\frac{1}{3}\o(N(Y,X),Z)\\
&\phantom{\nabla^o_X\o(Y,Z)cci}-\frac{1}{3}\o(Y,N(X,Z))-\frac{1}{3}\o(Y,N(Z,X))\\
&=\nabla^o_X\o(Y,Z)-\frac{1}{3}\o(N(X,Y),Z)-\frac{1}{3}\o(N(Y,X),Z)\\
&\phantom{\nabla^o_X\o(Y,Z)cci}-\frac{1}{3}\o(N(X,Y),Z)-\frac{1}{3}\o(N(Z,Y),X)\\
&=\left(1-\frac{1}{3}-\frac{1}{3}\right)\omega(N(X,Y),Z)+\frac{1}{3}\o(N(X,Z),Y)=0\;.
\end{aligned}
$$
\vskip-0.79cm\par\noindent
\qed
\vskip0.3cm\par\noindent
Let $(M,\o)$ be an almost symplectic manifold with an action of a Lie group $G$ preserving $\o$. It is a classical result of I.\ Vaisman that if $(M,\o)$ has a $G$-invariant linear connection then it also admits a $G$-invariant linear connection preserving $\o$ (see e.g. \cite[Corollary 1.3]{GRS}).

If $d\o=0$, Vaisman's result can be strengthened as follows.
\bp
\label{prop:invariantFedosov}
Let $(M,\o)$ be a symplectic manifold with an action of a Lie group $G$ which preserves $\o$. Assume there exists a $G$-invariant linear connection on $M$. Then $M$ has a torsion free $G$-invariant connection preserving $\o$.
\ep
\pf
Let $\nabla$ be a $G$-invariant linear connection on $M$. Symmetrisation (in the sense of Gelfand et al.) of this connection yields another connection
$$
\nabla^o_XY=\frac{1}{2}(\nabla_XY+\nabla_YX+[X,Y])
$$
which is $G$-invariant and torsion free. The corresponding $\o$-preserving torsion free connection from Proposition \ref{prop:fedosovcon} is $G$-invariant.
\qed
Let $(M = G/K,\omega)$ be a homogeneous symplectic manifold, on which a Lie group $G$ (not necessarily compact) acts effectively. Assume for simplicity that $G$ and $K$ are connected. We recall that the homogeneous manifold $M=G/K$ is called {\it reductive} if the Lie algebra $\gg$ of $G$ may be decomposed into a vector space direct sum $\gg=\gk+\gm$ of the Lie algebra $\gk$ of $K$ and a $\gk$-invariant subspace $\gm$.
\bp
\label{prop:redFedosov}
Any  reductive    homogeneous  symplectic manifold   $(M = G/K,\o)$ admits  a torsion free $G$-invariant connection which preserves $\o$, i.e., it is a homogeneous Fedosov manifold.
\ep
\pf   Any     such manifold  admits  a  $G$-invariant linear connection (for instance, the canonical connection \cite{KoNo}) and Proposition \ref{prop:invariantFedosov} applies.
\qed
\subsection{Uniqueness of homogeneous Fedosov structures}
\label{sec:lastsecII}
Let $(M=G/K,\o)$ be a $2n$-dimensional homogeneous symplectic manifold. For any $g\in G$, we denote by $L_{g}$ the corresponding left action on $M$. Let
\begin{equation}
\label{eq:Sbundle}
\pi:O_{\o}(M)\to M
\end{equation}
be the symplectic frame bundle of $M$, i.e,
the $\mathrm{Sp}_{n}(\bR)$-bundle of linear frames $u=(e_i)$ of the
tangent spaces of $M$ which are adapted to the symplectic form.

Given a symplectic frame $u_o=(e_i)\in O_{\o}(M)$ at $o=eK$, we let
$$P=G\cdot u_o=\Big\{g\cdot u_o:=(L_g{}_* e_i)\mid g\in G\Big\}\subset O_{\o}(M)$$
be the $G$-orbit of $u_o$. One gets in this way a natural reduction $\pi:P\rightarrow M$ of \eqref{eq:Sbundle}, usually called the {\it homogeneous $H$-structure associated with $M=G/K$}.

If $K_1$ is the kernel of the linear isotropy representation
$$
j:K\to \mathrm{Sp}(V)\;,\qquad V=T_oM\;,\qquad k\mapsto L_{k}{}_*|_{o}\,,
$$
then the map $g\mapsto g\cdot u_o$ determines a natural bundle isomorphism  between
\begin{equation}
\pi:G/K_1\rightarrow M=G/K
\label{eq:bundle2}
\end{equation}
and the $H$-structure $\pi:P\rightarrow M$, in such a way that the structure group $H=j(K)\cong K/K_1$ (see e.g. \cite{San}). We assume a fixed choice of $u_o$ and tacitly use it to identify $\pi:P\to M$ with \eqref{eq:bundle2}.
\\

We shall now restrict to the $4$-dimensional case. The following result is a straightforward consequence of Corollary \ref{thm:2}.
\begin{Lem}
Let $(M = G/K,\o)$ be a $4$-dimensional homogeneous symplectic manifold and assume the isotropy $H=j(K)\subset \mathrm{Sp}(V)$ has finite type. Then $K_1$ is discrete.
\end{Lem}

Our result deals with the uniqueness of homogeneous Fedosov structures compatible with extra geometric data.
\bt
\label{thm:last}
Let $(M = G/K,\o)$ be a homogeneous symplectic $4$-manifold
with finite type isotropy $H\subset \mathrm{Sp}(V)$. Assume there exists
a torsion free $G$-invariant connection $\nabla$ which is, in addition, compatible with the associated homogeneous $H$-structure $\pi:P\to M$ (in particular, $\nabla$ preserves the symplectic form $\o$). Then $K_1=\{1\}$ and
\begin{enumerate}
	\item $M=G/K$ is reductive;
	\item the connection $\nabla$ is unique.
\end{enumerate}
\et
\pf
An affine transformation of $M$ is the identity transformation if it leaves one linear frame fixed. This says $K_1=\{1\}$ and that
$\pi:P\to M$ is the $H$-principal bundle $\pi:G\to G/H$. It is a well known fact that a $G$-invariant connection on $\pi:G\to G/H$ exists if and only if the homogeneous manifold $M=G/H$ is reductive (see, e.g, \cite[Corollary 1.4.6]{CS}). This proves (1).
 
An invariant connection on $\pi:G \to G/H$ is completely determined by the associated Nomizu map, an $H$-equivariant map
$L:\gm\to\gh$, and it is torsion free if
$
L(x)y-L(y)x-\pi_{\gm}[x,y]=0
$,
where $\pi_{\gm}$ is the projection of $\gg$ on $\gm$ relative to the decomposition $\gg=\gh+\gm$ (see e.g. \cite{KoNo}).
  The difference of two Nomizu maps is an element of $\gh^{(1)}$, which is zero by Corollary \ref{thm:2}.
\qed
\subsection{Applications to symplectic Lie groups}
\label{sec:lastsecIII}
We conclude with an application of the results of  \S \ref{sec:lastsecI} and \S \ref{sec:lastsecII} to {\it symplectic Lie groups}, that is, the homogeneous symplectic manifolds $(M=G/K,\o)$ with trivial stabilizer $K=\{1\}$. It is clear that the existence of $\o$ is equivalent to $\gg$ being a {\it symplectic Lie algebra}, i.e., a Lie algebra endowed with a nondegenerate scalar $2$-cocyle.

Let $(G,\o)$ be a symplectic Lie group, of any (even) dimension. It is known that $G$ admits an invariant flat and torsion free connection $\nabla^o$ \cite{MR}. At the Lie algebra level, it can be described by
\begin{equation}
\label{eq:definitionlsp}
\o(\nabla^o_xy,z)=-\o(y,[x,z])\;,
\end{equation}
where $x,y,z\in\gg$ (we identify $\gg$ with the Lie algebra of left-invariant vector fields on $G$). It is customary to introduce the product
on $\gg$ given by covariant differentiation
\begin{equation}
\label{eq:lsp}
xy:=\nabla^o_xy
\end{equation}
and recast the flatness and zero-torsion conditions as
\begin{equation}
\label{eq:lspaxioms}
\begin{aligned}
(xy)z-x(yz)&=(yx)z-y(xz)\;,\\
[x,y]&=xy-yx\;.
\end{aligned}
\end{equation}
A vector space $A$ endowed with a bilinear product $xy$ satisfying the first equation in \eqref{eq:lspaxioms}, or equivalently for which the associator is symmetric in the first two variables, is called a {\it left-symmetric algebra}. Any left-symmetric algebra $A$ is Lie-admissible, that is, the commutator of the left-symmetric product is a Lie bracket. If this Lie bracket coincides with the Lie bracket of a Lie algebra structure $\gg$ given a priori on $A$, the left-symmetric product is called {\it compatible} (with $\gg$).
The construction of  \cite{MR} may be then summarized as follows.
\bp
Any symplectic Lie algebra $\gg$ admits a compatible structure \eqref{eq:lsp} of a left-symmetric algebra.
\ep
For more details on left-symmetric algebras, we refer the reader to the survey article \cite{Bu} and the references therein.
\vskip0.2cm\par
We collect a number of trace identities useful for our purposes. We denote left (resp. right) multiplication by $x\in\gg$ in the left-symmetric algebra by $L_x:y\mapsto xy$ (resp. $R_x:y\mapsto yx$). We let $(e_i)$ be a fixed symplectic basis of $\gg$ and $(e^i)$ the
dual basis of $\gg$ (i.e., $\o(e_i,e^j)=\delta_i^j$).
\bp
\label{prop:identitautili}
Let $\gg$ be a symplectic Lie algebra, with compatible left-symmetric structure. Then for all $x,y,z\in\gg$ we have:
\begin{enumerate}
	\item $\o(xy,z)+\o(zy,x)=0$,
	\item $\mathfrak{S}_{xyz}\Big\{\o(xy,z)\Big\}=0$,
	\item $\tr(R_x\circ R_y)=\tr(R_{xy})=2\tr(L_{xy})$,
	\item $\tr(R_x\circ R_y)=2\tr(R_y\circ L_x)=2\tr(R_x\circ L_y)$.
	\end{enumerate}
where $\mathfrak{S}_{xyz}$ is the cyclic sum over $x,y,z$.
\ep
\pf
The first property is immediate and the second is equivalent to the closure of $\o$, since $\nabla^o$ is torsion free. We now use the Einstein summation convention and compute
\begin{equation}
\begin{aligned}
\tr(R_{x})&=\o(e_ix,e^i)=-\o(e^ie_i,x)-\o(xe^i,e_i)\\
&=\o(e_ie^i,x)+\o(xe_i,e^i)=-\o(x,e_ie^i)+\tr(L_{x})\\
&=-\frac{1}{2}\o(x,[e_i,e^i])+\tr(L_{x})=\frac{1}{2}\o(e_ix,e^i)+\tr(L_{x})\\
&=\frac{1}{2}\tr(R_{x})+\tr(L_{x})\;,
\end{aligned}
\end{equation}
whence $\tr(R_x)=2\tr(L_x)$, for all $x\in\gg$; similarly
\begin{equation}
\begin{aligned}
\tr(R_x\circ R_y)&=\o((e_ix)y,e^i)\\
&=\o(e_i(xy),e^i)+\o((xe_i)y,e^i)-\o(x(e_iy),e^i)\\
&=\tr(R_{xy})+\tr(R_y\circ L_x)-\tr(L_x\circ R_y)\\
&=\tr(R_{xy})\;.
\end{aligned}
\end{equation}
Finally:
\begin{equation}
\begin{aligned}
\tr(R_x\circ R_y)&=\o((e_ix)y,e^i)\\
&=-\o(e^i(e_ix),y)-\o(ye^i,e_ix)\\
&=\o(y(e_ix),e^i)+\o(e_ix,ye^i)\\
&=\tr(L_y\circ R_x)-\o((ye^i)x,e_i)\\
&=\tr(L_y\circ R_x)+\tr(R_x\circ L_y)\;,
\end{aligned}
\end{equation}
which readily implies our last claim.
\qed
The bilinear forms $\tr(L_{xy})$ and $\tr(R_{xy})$ on a left-symmetric algebra have been extensively studied in the context of convex homogeneous cones \cite{Vin}, respectively, in connection with complete left-symmetric algebras \cite{He, Se}. To state the main result of this section, we need one last preliminary fact. Left multiplication in the left-symmetric algebra yields a representation $L:\gg\to\ggl(\gg)$ of the Lie algebra $\gg$ and we call the corresponding $\gg$-equivariant symmetric bilinear form
$$
\k(x,y):=\tr(L_x\circ L_y)
$$
the {\it left trace form}. We are not aware of any simple identity which relates the left trace form with the Killing form of the associated Lie algebra or with (any of) the bilinear forms considered in Proposition \ref{prop:identitautili}. Nonetheless, the following holds.
\bp
\label{eq:nilpotency}
Let $\gg$ be a symplectic Lie algebra, with the compatible left-symmetric structure. If the left trace form $\k=0$ then $\gg$ is solvable and, conversely, if $\gg$ is nilpotent then $\k=0$.
\ep
\pf
Iterating \eqref{eq:definitionlsp} we have
$$
\o(y,\ad_x^k(z))=(-1)^k\o(L_x^k(y),z)\;,
$$
for all positive integers $k$, where $x,y,z\in\gg$. In particular an element $x\in\gg$ is $\ad$-nilpotent if and only if the operator $L_x:\gg\to\gg$ is nilpotent. Note also that the kernel of the left regular representation $L:\gg\to\ggl(\gg)$ coincides with the center $\xi$ of $\gg$.

If $\k=0$, then we may apply Cartan's criterion to the matrix Lie algebra $L(\gg)\cong\gg/\xi$ and infer that the derived algebra of $\gg/\xi$ is nilpotent, whence $\gg$ is solvable.
If $\gg$ is nilpotent then the Lie algebra $L(\gg)\subset \ggl(\gg)$ consists of nilpotent operators  and by Engel's theorem it can be represented by upper triangular matrices. Hence $\k=0$.
\qed

Let now $(G,\o)$ be a symplectic Lie group and $\nabla$
the torsion free symplectic connection associated to the left-symmetric product \eqref{eq:lsp} on $\gg$ via Proposition \ref{prop:fedosovcon}. The section $N$ of $\gg^*\otimes \gg^*\otimes \gg$ satisfies
$$
\begin{aligned}
\omega(N(x,y),z)&=\nabla_x^o\o(y,z)\\
&=x(\o(y,z))-\o(xy,z)-\o(y,xz)\\
&=-\o(xy,z)+\o([x,y],z)\\
&=-\o(yx,z)\;,
\end{aligned}
$$
whence $N(x,y)=-yx$ and
$$
\begin{aligned}
\nabla_xy&=\nabla^o_xy+\frac{1}{3}N(x,y)+\frac{1}{3}N(y,x)\\
&=\frac{2}{3}xy-\frac{1}{3}yx\;,
\end{aligned}
$$
for all $x,y\in\gg$.
\bp
The curvature of $\nabla$ is given by
$$
R(x,y)=-\frac{1}{9}[R_x,R_y]-\frac{2}{9}L_{[x,y]}+\frac{1}{9}R_{[x,y]}\;,
$$
for all $x,y\in\gg$.
\ep
\pf
We compute
$$
\begin{aligned}
R(x,y)z&=\nabla_x\nabla_y z-\nabla_y\nabla_xz-\nabla_{[x,y]}z\\
&=\nabla_x(\frac{2}{3}yz-\frac{1}{3}zy)-\nabla_y(\frac{2}{3}xz-\frac{1}{3}zx)-\nabla_{xy}z+\nabla_{yx}z\\
&=\frac{4}{9}(x(yz)-y(xz))
+\frac{2}{9}((xz)y-x(zy))
-\frac{2}{9}((yz)x-y(zx))
-\frac{1}{9}(zx)y
+\frac{1}{9}(zy)x
\\
&\;\;\;\,-\frac{2}{3}(xy)z+\frac{2}{3}(yx)z+\frac{1}{3}z(xy)-\frac{1}{3}z(yx)\\
&=\frac{1}{9}(zx)y-\frac{1}{9}(zy)x-\frac{2}{9}[x,y]z+\frac{1}{9}z[x,y]\;,
\end{aligned}
$$
where in the last step we repeatedly used the left-symmetric property of the product.
\qed
We now state the main result of this section.
\bt
\label{thm:nilpoliegroup}
Let $(G,\o)$ be a symplectic Lie group (of any dimension) with Lie algebra $\gg$ and $$\nabla_xy=\frac{2}{3}xy-\frac{1}{3}yx\;,$$
$x,y\in\gg$, the  symplectic torsion free connection naturally associated to the compatible structure of left-symmetric algebra on $\gg$ (cf. \S \ref{sec:lastsecI}). Then the Ricci tensor of $\nabla$ is given by
$$
\operatorname{ric}^{\nabla}(x,y)=\frac{1}{9}\Big(\tr(L_{xy})+\tr(L_x\circ L_y)\Big)\;.
$$
If $(G,\o)$ is a nilpotent symplectic Lie group then $\operatorname{ric}^{\nabla}=0$.
\et
\pf
We fix a symplectic basis $(e_i)$ of $\gg$ with dual basis $(e^i)$, use the Einstein summation convention and compute
\begin{equation}
\begin{aligned}
9\operatorname{ric}^{\nabla}(x,y)&=9\,\tr\Big\{z\mapsto R(x,z)y\Big\}\\
&=9\,\o(R(x,e_i)y,e^i)\\
&=+\o((yx)e_i,e^i)-\o((ye_i)x,e^i)
-2\o((xe_i)y,e^i)\\
&\;\;\;\,+2\o((e_ix)y,e^i)+\o(y(xe_i),e^i)
-\o(y(e_ix),e^i)\\
&=\tr(L_{yx})-\tr(R_x\circ L_y)-2\tr(R_y\circ L_x)\\
&\;\;\;\,+2\tr(R_y\circ R_x)+\tr(L_y\circ L_x)-\tr(L_y\circ R_x)
\;.
\end{aligned}
\end{equation}
Using Proposition \ref{prop:identitautili}, one readily gets the formula for the Ricci tensor. If $G$ is nilpotent, then each $L_x:\gg\to\gg$ is nilpotent too (see the proof of Proposition \ref{eq:nilpotency}) and $\operatorname{ric}^\nabla=0$ follows from $\tr(L_{xy})=0$ and Proposition \ref{eq:nilpotency}.
\qed
\section*{Acknowledgments}
\vskip0.2cm\par\noindent
D.\ A.\ acknowledges support of the grant n. 18-00496S of
the Czech Science Foundation. The research of A.\ S.\ is supported by the
project ``Lie superalgebra theory and its applications'' of the University of Bologna and partly supported
by the Project Prin 2015 ``Moduli spaces and Lie Theory''.

\vskip0.2cm\par\noindent
Dmitri V.\ Alekseevsky\\
A. A. Kharkevich Institute for Information Transmission Problems, B. Karetnui per. 19, 127051, Moscow, Russia and
University of Hradec Kr\'alov\'e,
Faculty of Science, Rokitansk\'eho 62, 50003 Hradec Kr\'alov\'e,  Czech
Republic. email: dalekseevsky@iitp.ru
\vskip0.3cm\par\noindent
Andrea Santi\\
Dipartimento di Matematica, Università di Bologna, Piazza di Porta San Donato 5, 40126, Bologna, Italy.
email: asanti.math@gmail.com
\enddocument